\documentclass[reqno,11pt]{amsart}
\usepackage{color}

\newtheorem{theorem}{Theorem}[section]
\newtheorem{lemma}[theorem]{Lemma}
\newtheorem{corollary}[theorem]{Corollary}

\theoremstyle{definition}
\newtheorem{assumption}[theorem]{Assumption}
\newtheorem{definition}[theorem]{Definition}

\theoremstyle{remark}
\newtheorem{remark}[theorem]{Remark}

\makeatletter
\def\dashint{\operatorname%
{\,\,\text{\bf--}\kern-.98em\DOTSI\intop\ilimits@\!\!}}
\makeatother

\newcommand\bC{\mathbb{C}}

\newcommand\bL{\mathbb{L}}
\newcommand\bM{\mathbb{M}}
\newcommand\bR{\mathbb{R}}

\newcommand\bS{\mathbb{S}}
\newcommand\bZ{\mathbb{Z}}

\newcommand\cF{\mathcal{F}}

\newcommand\cM{\mathcal{M}}

\newcommand\sfu{{\sf u}}

\newcommand{\loc}{{\rm loc}}

 \newcommand{\mysection}[1]{\section{#1}
 \setcounter{equation}{0}}

\begin{document}

\title[Fully nonlinear elliptic and parabolic equations]
{Fully nonlinear elliptic and parabolic equations in weighted and mixed-norm Sobolev spaces}
\author[H. Dong]{Hongjie Dong}
\address[H. Dong]{Division of Applied Mathematics, Brown University,
182 George Street, Providence, RI 02912, USA}
\email{Hongjie\_Dong@brown.edu}
\thanks{H. Dong was partially supported by the NSF under agreement DMS-1600593.}
\author[N. V. Krylov]{N.V. Krylov}
\address[N. V. Krylov]{127 Vincent Hall, University of Minnesota,
 Minneapolis, MN, 55455}
\email{nkrylov@umn.edu}

\begin{abstract}
We prove weighted and mixed-norm Sobolev estimates for fully nonlinear elliptic and parabolic equations in the whole space under a relaxed convexity condition with almost VMO dependence on space-time variables. The corresponding interior and boundary estimates are also obtained.
\end{abstract}

\maketitle

\mysection{Introduction}

The goal of this paper is to establish weighted and mixed-norm Sobolev estimates
for fully nonlinear second-order elliptic and parabolic equations with almost VMO
dependence on space-time variables,  under a relaxed convexity condition. The interest in results concerning equations in spaces with mixed Sobolev norms arises, for example, when one wants to have better regularity of traces of solutions of parabolic equations
for each time slice while treating linear or nonlinear equations.

The usual Sobolev space theories of
{\em linear\/} elliptic and parabolic equations
with continuous main coefficients
has long and rich history reflected in
lots of papers and books. In early nineties
Chiarenza, Frasca, and Longo, and Bramanti and
Cerutti discovered a way which allows
main coefficients to be almost in VMO
rather  than continuous. Their approach was also continued in quite a few papers and books.
As the previous theory, this approach is based
on the theory of singular integrals or
its versions and explicit integral representation
of  solutions of model equations. The same approach
  also  works for equations with sufficiently regular
coefficients in Sobolev spaces with  Muckenhoupt $A_{p}$-weights,
as is shown, for instance, in \cite{CD} and the
references therein. About ten years ago a different approach
was suggested based on the Fefferman-Stein theorem
in place of the theory of singular integrals.
This approach is more flexible and applies
to nonlinear equations as well as to
linear ones and does not require any explicit
representation of solutions in any model
case. For instance, it allowed the authors of
\cite{DK16, DG18} to generalize the results of the type in
\cite{CD} to a   large extent to a very wide
range of equations with almost VMO
coefficients and, in addition, also derive
mixed norms estimates with $A_{p}$-weights.

Our goal is to prove similar results for fully nonlinear equations.

Let $\bR^{d}$ be the $d$-dimensional Euclidean
space of points $x=(x_{1},\ldots,x_{d})$ and $\bS$ be the set
of $d\times d$ symmetric matrices.
For $\delta\in(0,1)$, by $\bS_{\delta}$ we denote the subset
of $\bS$ consisting of matrices whose eigenvalues are
between $\delta$ and $\delta^{-1}$. We are interested in elliptic operators in the form
$$
F[u]:=F(D^2 u,x),
$$
where $F=F(\sfu'',x),\sfu''\in \bS,x\in \bR^d$, is a given function,
as well as the corresponding parabolic operators in the form
$$
\partial_t u+F[u]:=\partial_t u+F(D^2 u,t,x).
$$
Here  and everywhere below
$$
D^2u = (D_{ij}u),\quad Du = (D_iu),\quad D_i =\frac{\partial}
{\partial x_i},\quad D_{ij} = D_iD_j,\quad \partial_t =\frac{\partial}
{\partial t}.
$$

Under the assumption that $F$ is Lipschitz continuous with respect to $\sfu''$, $F(0,x)=0$, $F$ is almost convex in $\sfu''$ and almost VMO in $x$ for large values of $|\sfu''|$, we obtain weighted Sobolev estimates in the whole space with Muckenhoupt $A_p$-weights. See Section \ref{section 4.10.3} and Theorem \ref{theorem 12.15.1} for more precise assumptions and the result. By using a powerful extrapolation theorem due to J. L. Rubio de Francia \cite{MR745140}, we then derive mixed-norm Sobolev estimates in the whole space under some additional conditions. See Theorem \ref{theorem 4.9.2}. For operators $F$ which are positive homogeneous of degree one with respect to $\sfu''$, we prove a local mixed-norm estimate. See Theorem \ref{theorem 12.29.1}. We also consider fully nonlinear elliptic equations in half spaces, and prove estimates near the boundary with $A_{p}$-weights on $\bR^d_+$ and, as a typical example, weights which are powers
of the distance to the boundary. See Sections \ref{section 4.10.1} and \ref{section 4.10.2}. The corresponding estimates for parabolic equations in the whole space, half spaces, balls, and half balls are also established in Sections \ref{section 4.23.1} and \ref{sec7}. It is worth noting that one can also consider operators $F$ with lower order terms. However, in order not to overburden this paper, we only consider operators which depends only on $D^2 u$ and $x$ (and also $t$ in the parabolic case).

Our proofs of weighted estimates are based on mean oscillation estimates proved earlier in \cite{Kr_13,Kr_18}, the Hardy-Littlewood maximal function theorem, and a local version of the Fefferman-Stein sharp functions theorem with $A_p$-weights, which is one of our main results and is stated in Corollary \ref{corollary 12.3.1} below. Such local version of the Fefferman-Stein sharp functions theorem allows us to derive estimates without relying on a partition of unity argument, which is not applicable to general fully nonlinear operators. The key ingredients in the proof of mean oscillation estimates in \cite{Kr_13,Kr_18} are the Evans-Krylov theorem  and a $W^2_\varepsilon$
estimate for equations with measurable coefficients, which is originally due to F.H. Lin \cite{Lin86}.
For mixed-norm estimates, we follow the argument in \cite{DK16} by using a generalized extrapolation theorem, Theorem \ref{theorem 12.2.1}, in the spirit of J. L. Rubio de Francia \cite{MR745140}.

The interior (usual)
$W^2_p$ estimates for fully nonlinear elliptic equations  were derived in \cite{CC_95}, basically,
under the convexity
assumption on $F$ with respect to $\sfu''$
and almost continuity assumption with respect to $x$.
In \cite{Wi09} global estimates  were obtained under the same kind
of assumption. These results  were obtained by using the theory of viscosity solutions. The same theory applied in \cite{CKS00} to parabolic case yields
similar results under similar assumptions as in the elliptic case.

For elliptic Bellman's
equations with {\em VMO dependence\/}
on the independent variables
the interior
$W^2_p$ estimates were first obtained in \cite{Kr_10}.

Later, boundary and similar estimates for parabolic equations, as well as a solvability result, were obtained in \cite{DKL_12}. The {\em
relaxed convexity\/} and VMO conditions (Assumptions \ref{assump1} and \ref{assump1b}) in the current paper are adopted from \cite{Kr_13}, in which the existence of $W^2_p$ solutions for fully nonlinear elliptic equations in domains was proved. See also \cite{Kr_18} for a result for parabolic equations.
This paper is a continuation of this line of research in the weighted and mixed-norm settings. For other relevant results in the literature, we refer the reader to \cite{Kr_10, DKL_12, Kr_13, Kr_18} and a recent book \cite{Kr17} by the second named author.

The remaining part of the paper is organized as follows. In the next section, we recall some definitions and facts from Chapter 3 of \cite{Kr_08} and prove a local version of the Fefferman-Stein sharp function theorem. We consider elliptic equations in the whole space and in balls in Section \ref{section 4.10.3}, and in half spaces and in half balls in Sections \ref{section 4.10.1} and \ref{section 4.10.2}. In Sections \ref{section 4.23.1} and \ref{sec7}, we prove analogous results for parabolic equations. In the appendix, we state and prove a generalized extrapolation theorem, Theorem \ref{theorem 12.2.1}.

\mysection{Partitions and sharp functions}

 For reader's convenience, we first recall some definitions and facts from
Chapter 3 of \cite{Kr_08}.
Let $(\Omega,\cF,\mu)$ be a complete measure space with
a $\sigma$-finite measure $\mu$, such that
$$
\mu(\Omega)=\infty.
$$
Let $\cF^{0}$ be the subset
 of $\cF$ consisting of all
sets $A$ such that $\mu(A)<\infty$.
By $\bL$ we denote a fixed dense subset
of $L_{1}(\Omega)=L_{1}(\Omega,\cF,\mu)$.
For any $A\in\cF$ we
set
$$
|A|=\mu(A).
$$
For $A\in\cF^{0}$ and   functions $f$
summable on $A$ we
  use the notation
$$
f_{A}=\dashint_{A}f\,\mu(dx):=\frac{1}{|A|}\int_{A}f(x)\,\mu(dx)
\quad \bigg(\frac{0}{0}:=0\bigg)
$$
 for  the average value of $f$ over $A$.
We write
$f\in L_{1,\loc}(\Omega)$ if $fI_{A}\in L_{1}(\Omega)$
for any $A\in\cF^{0}$.

\begin{definition}
                                          \label{definition 5.30.1}
Let $\bZ=\{n: n=0,\pm1,\pm2,\ldots\}$ and let
 $(\bC_{n},n\in\bZ)$ be
a sequence of   partitions of $\Omega$ each consisting of
countably many disjoint
  sets $C\in\bC_{n} $ and such that $\bC_{n}\subset
\cF^{0}$ for each $n$.
For each $x\in\Omega$ and $n\in \bZ$ there exists (a unique)
$C\in\bC_{n}$ such that $x\in C$. We denote
this $C$ by
$C_{n}(x)$.

We call the sequence $(\bC_{n},n\in\bZ)$ {\em a filtration
of partitions\/}
if the following conditions are satisfied.

(i) The elements of partitions
are ``large"     for big negative $n$'s and ``small"
for  big positive
$n$'s:
$$
\inf_{C\in\bC_{n}}|C|\to\infty\quad\text{as}\quad n\to-\infty,\quad
 \lim_{n\to\infty}f_{C_{n}(x)}=f(x)\quad\text{(a.e.)}\quad
\forall f\in\bL.
$$

(ii) The partitions are nested: for each $n$ and $C\in\bC_{n}$
there is a (unique) $C'\in\bC_{n-1}$
 such that $C\subset C'$.

(iii) The following regularity property holds: for any $n$, $C$, and $C'$
as in (ii) we have
$$
|C'|\leq N_{0}|C|,
$$
 where $N_{0}$
is a constant independent of $n,C, C'$.
\end{definition}

We set
$$
\bC_{\infty}=\bigcup_{n}\bC_{n}.
$$

\begin{definition}
                                          \label{definition 5.30.2}
 Let $\bC_{n}$, $n\in\bZ$, be a filtration of partitions
of $\Omega$.

(i) Let $\tau=\tau(x)$ be a function on $\Omega$ with values in
$\{\infty,0,\pm1,\pm2,\ldots\}$. We call $\tau$ {\em a stopping time\/}
(relative to the filtration)
if, for each $n=0,\pm1,\pm2,\ldots$, the set
$$
\{x:\tau(x)=n\}
$$ is either empty or else is the union of
some elements of $\bC_{n}$.

(ii) For a function $f\in L_{1,\loc}(\Omega)$ and
$n\in\bZ$,
we denote
$$
f_{|n}(x)= \dashint_{C_{n}(x)}f(y)\,\mu(dy).
$$

If we are also given a stopping time $\tau$, we
let
$$
f_{|\tau}(x) =f_{|\tau(x)}(x)
$$
 for those $x$ for which
$\tau(x)<\infty$  and $f_{|\tau}(x)=f(x)$ otherwise.
\end{definition}

The simplest example of a stopping time is given by
$\tau(x)\equiv0$.

We are going to use the following simple properties of the objects
introduced above.

\begin{lemma}
                                 \label{lemma 06.5.30.1}
 Let $\bC_{n}$, $n\in\bZ$, be a filtration of partitions
of $\Omega$.

(i) Let   $f\in L_{1,\loc}(\Omega)$, $f\geq0$,
 and let $\tau$ be a stopping time.
  Then
\begin{gather}
                                          \label{06.5.1.1}
\int_{\Omega} f_{|\tau}(x) I_{\tau<\infty}\,\mu(dx)=
\int_{\Omega} f(x)  I_{\tau<\infty}\,\mu(dx),
\\[10pt]
                                          \label{07.9.27.2}
\int_{\Omega} f_{|\tau}(x)  \,\mu(dx)=
\int_{\Omega} f(x)   \,\mu(dx).
\end{gather}

(ii) Let $g\in L_{1}(\Omega)$, $g\geq0$,
and let $\lambda>0$ be a constant. Then
\begin{equation}
                                          \label{06.5.1.2}
\tau(x):=\inf\{n:g_{|n}(x)>\lambda\}\quad(\inf\emptyset:=\infty)
\end{equation}\vskip .1in \noindent
is a stopping time. Furthermore,
we have
\begin{equation}
                                          \label{06.5.1.3}
0\leq g_{|\tau}(x)I_{\tau<\infty}\leq N_{0}\lambda,\quad
|\{x:\tau(x)<\infty\}|\leq  \lambda^{-1}\int_{\Omega}g(x)
I_{\tau<\infty}\,\mu(dx).
\end{equation}

\end{lemma}

 Define
{\em the maximal function\/}
of $f$ by
$$
\cM f(x)=\sup_{n<\infty}|f|_{|n}(x),
$$
 so that
$\cM f=\cM |f|$.

Notice that Lemma \ref{lemma 06.5.30.1}
implies the following.

\begin{corollary}[Maximal inequality]
                                          \label{corollary 06.5.31.8}
For   $\lambda>0$
and   nonnegative $g\in L_{1}(\Omega)$, the
{\em maximal inequality\/}
 holds:
\begin{equation}
                                                 \label{06.06.29.2}
|\{x:\cM g(x)>\lambda\}|\leq \lambda^{-1}\int_{\Omega}
g(x)I_{\cM g>\lambda}\,\mu(dx).
\end{equation}\vskip .1in \noindent
\end{corollary}

Indeed, for $\tau$ as in \eqref{06.5.1.2}, we have
$$
\{x:\cM g(x)>\lambda\}=\{x:\tau(x)<\infty\}.
$$

\begin{corollary}
                                          \label{corollary 06.5.31.9}
Let $p\in(1,\infty)$, $g\in L_{1}(\Omega)$,
$g\geq0$. Then
$$
 \|\cM g\|_{L_{p}(\Omega)}\leq q\|g\|_{L_{p}(\Omega)},
$$
where $q=p/(p-1)$.
\end{corollary}

The following extends Corollary \ref{corollary 06.5.31.9} to
$g\in L_{p}(\Omega)$ .

\begin{theorem}
                                         \label{theorem 06.06.29.1}
For any
$p\in(1,\infty)$ and $g\in L_{p}(\Omega)$,
$$
 \|\cM g\|_{L_{p} (\Omega)}\leq q\|g\|_{L_{p}(\Omega)}.
$$
\end{theorem}

Let $w=w(x)$ be a nonnegative function on $\Omega$,
such that $\chi(C)<\infty$ for any $C\in \bC_{\infty}$, where
$$
\chi(A):=\int_{A}w\,\mu(dx) .
$$
For $ \beta\in(0,1]$,  we say that $w$ is of $\beta$-type if
$$
\frac{\chi(A)}{\chi(C)}\leq N_{w,\beta} \frac{|A|^{\beta}}{|C|^{\beta}}
$$
for any measurable $A\subset C$ and $C\in \bC_{\infty}$,
where $N_{w,\beta}$ is a (finite) constant independent of $C$ and $A$.
\begin{remark}
                                         \label{remark 5.27.1}

In some of our applications $\Omega$ will be a linear metric space
with filtration of either dyadic standard   or parabolic cubes
and $w$ will be an $A_{p}$-weight with respect to the
corresponding metric. One knows that in such situations
if $w\in A_{p}$ and $[w]_{p}\leq K_{0}$,
where $K_{0}$ is a constant, then $w$
is of $\beta$-type for an appropriate $\beta$
 and $N_{w,\beta}$ both depending only on $K_{0}$ and the metric.
\end{remark}

The following is a combination of Theorem 2.5 of \cite{DK16} and
Lemma 5.1 of \cite{Kr_10}.

\begin{lemma}
                                         \label{lemma 12.16.1}
Let $\gamma\in(0,1]$, $v\in L_{1,\loc}(\Omega)$,
and let $v_{|n}\to0$ as $n\to-\infty$ on $\Omega$. Assume that
  $|u|\leq v$
and
for any   $C\in\bC_{\infty}$ there exists a
measurable function  $u^{C}$ given on $C$ such that
$|u|\leq u^{C}\leq v$ on $C$ and, for any $x\in C$
\begin{equation}
                                                      \label{6.29.3}
\Big( \dashint_{C}
\dashint_{C}\big|u^{C}(z)-u^{C}(y)\big|^{\gamma}\,\mu(dz)\mu(dy)\Big)^{1/\gamma} \leq
 g (x)\ .
\end{equation}
Let $w$ be of $\beta$-type. Then for any $\lambda>0$ we have
\begin{equation}
                                                      \label{6.29.1}
 \chi\big\{x:\big|u(x)\big|\geq\lambda\big\} \leq N_{w,\beta}
\nu^{-\beta}\lambda^{-\gamma\beta}
\int_{\Omega}g^{\gamma\beta}(x)I_{\cM v(x)>\alpha \lambda}\,\chi(dx),
\end{equation}
where $\alpha=(2N_{0})^{-1}$ and $\nu= 1-2^{-\gamma}$.
\end{lemma}

Proof. Obviously we may assume that   $u\geq0$.
Fix a $\lambda>0$ and define  \vspace{5pt}
$$
\tau(x)=\inf\big\{n\in\bZ:v_{|n}(x)>\alpha\lambda\big\}.
\vspace{5pt}$$
We know that $\tau$ is a stopping time and if $\tau(x)<\infty$,
then
$$
  v_{|n}(x)\leq \lambda/2,\quad\forall n\leq\tau(x).
$$
We also know that $v_{|n}\to v\geq u$ (a.e.) as $n\to\infty$
(the Lebesgue differentiation theorem).
It follows  that (a.e.)
$$
\big\{x:u(x)\geq\lambda\big\}=\big\{x:u(x)\geq\lambda,\tau(x)<\infty\big\}
\vspace{5pt}$$
$$
=\big\{x:u(x)\geq\lambda,  v_{|\tau}(x)\leq \lambda/2\big\}
=\bigcup_{n\in\bZ}\bigcup_{C\in
\cF^{\tau}_{n}}A_{n}(C),
$$
where
$$
A_{n}(C):=\big\{x\in C:u(x)\geq\lambda, v_{|n}(x)\leq \lambda/2\big\},
\vspace{5pt}$$
and $\cF^{\tau}_{n}$ is the family of disjoint elements
of $\bC_{n}$ such that
$$
\big\{x:\tau(x)=n\big\}=\bigcup_{C\in \cF^{\tau}_{n}}C.
$$

Next, for each $n\in\bZ$ and $C\in\bC_{n}$ on the
set $A_{n}(C)$,
if it is not empty,
we have $v_{|n}=v_{C}$ and on $A_{n}(C)$
$$
u^{\gamma}-(v_{C})^{\gamma}\geq\lambda^{\gamma}(1-2^{-\gamma})
=\nu\lambda^{\gamma} .
 $$
  We use this
and the inequality $|a-b|^{\gamma}\geq|a|^{\gamma}-|b|^{\gamma}$
and conclude that for $x\in A_{n}(C)$,
$$
\dashint_{C}\big|u^{C}(x)-u^{C}(y)\big|^{\gamma}\,\mu(dy)
\geq\big(u^{C}(x)\big)^{\gamma}-\dashint_{C}\big(u^{C}(y)\big)^{\gamma}\,\mu(dy)
 $$
$$
\geq u^{\gamma}(x)-\dashint_{C}v^{\gamma}(y)\,\mu(dy)
\geq u^{\gamma}(x)-\big(v_{C}(x)\big)^{\gamma}
\geq\nu\lambda^{\gamma},
\vspace{5pt}$$
 so that
by Chebyshev's inequality  \vspace{5pt}
$$
\big|A_{n}(C)\big|\leq
\nu^{-1}\lambda^{-\gamma}\int_{ C }
\dashint_{C}\big|u^{C}(z)-u^{C}(y)\big|^{\gamma}\,\mu(dz)\mu(dx).
\vspace{5pt}$$

It follows by assumption \eqref{6.29.3} that \vspace{5pt}
$$
\frac{\big|A_{n}(C)\big|}{|C|}\leq
\nu^{-1}\lambda^{-\gamma}g^{\gamma}(x)
\vspace{5pt}$$
for any $x\in\Omega$.  Since $w$ is of $\beta$-type,
$$
\chi(A_{n}(C))\leq N_{w,\beta}\nu^{-\beta}\lambda^{-\gamma\beta}
g^{\gamma\beta}(x)\chi(C).
$$
Since this holds for any $x\in C$,
$$
\chi(A_{n}(C))\leq N_{w,\beta}\nu^{-\beta}\lambda^{-\gamma\beta}\int_{C}
g^{\gamma\beta}(x)\chi(dx).
$$
Hence,
$$
\chi\big\{x:u(x)\geq\lambda\big\}\leq N_{w,\beta}\nu^{-\beta}\lambda^{-\gamma\beta}
\sum_{n\in\bZ}\sum_{C\in \cF^{\tau}_{n}}\int_{C}g ^{\gamma\beta}\,\chi(dx)
\vspace{5pt}$$
$$
=N_{w,\beta}\nu^{-\beta}\lambda^{-\gamma\beta}\int_{\Omega}g^{\gamma\beta}I_{\tau<\infty}\,
\chi(dx).
\vspace{5pt}$$
It only remains to observe that $\{\tau<\infty\}=\{\cM
v>\alpha\lambda\}$. The lemma is proved.  \qed

\begin{corollary}
                                                 \label{corollary 12.14.1}
Under the assumption of  Lemma \ref{lemma 12.16.1}, for any $p>\gamma \beta$,
$$
\int_{\Omega}|u|^{p} \,\chi(dx)
\leq N\Big(\int_{\Omega}|\cM v|^{p}\,\chi(dx)\Big)^{(p-\gamma\beta)/p}
\Big(\int_{\Omega}|g|^{p}\,\chi(dx)\Big)^{\gamma\beta/p},
$$
where $N$ depends only on
$N_0$, $N_{w,\beta}$, $p$, $\beta$, and $\gamma$.
\end{corollary}
 Indeed, by Lemma \ref{lemma 12.16.1} and the Fubini theorem,
$$
\int_{\Omega}|u|^{p}\,\chi(dx)=
p\int_0^\infty \chi\big\{x :|u(x)|\geq\lambda\big\}\lambda^{p-1}\,d\lambda
$$
$$
\le pN_{w,\beta}\nu^{-\beta}\int_0^\infty \int_\Omega g^{\gamma\beta}(x)I_{\cM v(x)>\alpha \lambda}\lambda^{p-1-\gamma \beta}\chi(dx)\,d\lambda
$$
$$
=pN_{w,\beta}\nu^{-\beta}/(p-\gamma \beta)\int_\Omega g^{\gamma\beta}(x) (\cM v(x)/\alpha)^{p-\gamma\beta}\,\chi(dx).
$$
To get the desired inequality, it only remains
to apply H\"older's inequality.

For $m\in \bZ$ introduce
$$
u^{\#}_{\gamma,m}(x)=\sup_{n\geq m}
\sup_{\substack{C\in \bC_{n},\\
C\ni x}} \Big(\dashint_{C}
 \dashint_{C}|u(z)-u(y)|^{\gamma}\mu(dz)\mu(dy)\Big)^{1/\gamma},
$$
$$
\cM_{m}v=\sup_{n\leq m}v_{\mid n}.
$$

\begin{corollary}
                                           \label{corollary 12.3.1}
Take $m\in \bZ$. Assume that $|u|_{\mid n}\to0$
as $n\to-\infty$, and let $w$ be of $\beta$-type.  Then for any
$p>\gamma\beta$,
$$
\int_{\Omega}|u|^{p}\,\chi(dx)\leq NI^{(p-\gamma\beta)/p}
J^{\gamma\beta/p},
$$
where
$$
I=\int_{\Omega}|\cM u|^{p}\,\chi(dx),
$$
$$
J=\int_{\Omega}\big(  u^{\#}_{\gamma,m}
+\cM^{
1/\gamma}_{m}(|u|^{\gamma})\big)^{p}\,\chi(dx),
$$
and the constant $N$ depends only on
$N_0$, $N_{w,\beta}$, $p$, $\beta$, and $\gamma$.
\end{corollary}

This obviously follows from Corollary \ref{corollary 12.14.1}
with $u^C=v=|u|$
since for $n\leq m$ the left-hand side of \eqref{6.29.3}
is less that $2^{1/\gamma}\cM^{1/\gamma}_{m}(v^{\gamma})$.

\mysection{Elliptic case}
                                   \label{section 4.10.3}

In this Section, we study fully nonlinear elliptic equations in weighted and mixed-norm Sobolev spaces.
Set
$$
B_{r}(x)=\{y\in\bR^{d}:|x-y|<r\},\quad B_{r}=B_{r}(0).
$$
Suppose that we are given a function $F(\sfu'',x)$, $\sfu''\in\bS$, $x\in\bR^{d}$. In our results we will impose some of the following
assumptions.

\begin{assumption}[$\theta$]
                                \label{assump1}
  (i) The function $F$ is Lipschitz continuous with respect to $\sfu''$
 with Lipschitz constant $K_{F}$ and  $F(0,x)\equiv0$.

There exist   $R_0\in(0,1]$  and $\tau_{0}\in[0,\infty)$
such that, if   $r\in (0, R_0]$ and $z\in \bR^{d}$, then
one can find a {\em convex\/} function $\bar{F} (\sfu'' )=
\bar{F}_{z,r } (\sfu'' )$ (independent
of $x $)  for which

(ii) We have $\bar{F}(0)=0$ and
$ D_{\sfu'' }\bar{F} \in\bS_{\delta}$
at all points of differentiability
of $\bar{F}$;

(iii)
For any $\sfu''\in\bS$ with $|\sfu''|=1$,
 we have
\begin{equation}
                                                \label{7.30.2}
\int_{ B_{r}(z)}\sup_{\tau>\tau_{0}}\tau^{-1}
\big|F\big( \tau \sfu'' ,x\big)-\bar{F}(\tau \sfu'')\big| \,dx\leq
\theta
\big| B_{r}(z)\big|,
\end{equation}
where by $|A|$ we denote the volume of $A$ in $\bR^{d}$.
\end{assumption}

\begin{assumption}
                                         \label{assumption 4.26.1}

The function $F$ is Lipschitz continuous with respect to $\sfu''$,
$F(0,x)\equiv 0$,
and $ D_{\sfu'' }F \in\bS_{\delta}$
at all points of differentiability
of $F$.
\end{assumption}

\begin{remark}
                                          \label{remark 4.26.1}
Assumption \ref{assumption 4.26.1} implies that,
for any $ \sfu''\in\bS$ and $x\in\bR^{d}$, we have $F( \sfu'',x)
=a^{ij} \sfu''_{ij}$, where $a=(a^{ij})\in\bS_{\delta}$.
\end{remark}

For functions $h$ on $\bR^{d}$, $\rho>0$, and $x\in
\bR^{d}$, introduce
$$
h^{\sharp}_{ \gamma,\rho}(x) =
\sup_{\substack{r\in(0,\rho],\\
B_{r}(x_{0})\ni x}}\Big(\dashint_{B_{r}(x_{0}) }
\dashint_{B_{r}(x_{0}) }
\big|h(x_{1})-h(x_{2})\big|^{\gamma}\,dx_{1}dx_{2}\Big)^{1/\gamma},
$$
\begin{equation}
                                                    \label{7,2,1}
\begin{split}
\bM  h(x) &=\sup_{\substack{r>0,\\
B_{r}(x_{0})\ni x}}\dashint_{B_{r}(x_{0}) }
|h(y)|\,dy,\\
\bM_{\rho}  h(x) &=\sup_{\substack{r\in [\rho,\infty) ,\\
B_{r}(x_{0})\ni x}}\dashint_{B_{r}(x_{0}) }
|h(y)|\,dy.
\end{split}
\end{equation}

We set $\Omega=\bR^{d}$ and
for $n\in\bZ$ we take $\bC_{n}$
as the collection of $x+[0,2^{-n})^{d}$, $x\in 2^{-n}\bZ^{d}$.
We also set $\mu$ to be Lebesgue measure
and $\bL$ to be the set of continuous functions with compact support.
Then
observe that  for a constant $c=\sqrt d/2$,
\begin{equation}
                                            \label{12.15.5}
h^{\#}_{ \gamma,m}\leq Nh^{\sharp}_{ \gamma,c2^{-m}},
\quad \cM_{m}h\leq  N \bM_{c2^{-m}}h.
\end{equation}

From Lemma 3.6 of \cite{Kr_13} and the proof of
Lemma 5.2 (related to estimates in bounded domains)
of \cite{Kr_13} one can easily obtain the following result.

\begin{lemma}
                                \label{lemma 12.15.1}
Let $u\in W^{2}_{d,\loc}(\bR^{d})$, $\mu\in(0,\infty)$, $\nu\geq2$,
$\xi\in(1,\infty)$.  Then there exists
$\theta=\theta(d,\delta,K_{F},\mu,\xi )\in(0,1)$
such that, if Assumption \ref{assump1} $(\theta)$
is satisfied, then one can find $
 \gamma_0=\gamma_0(d,\delta)
\in(0,1)$, $ \alpha=\alpha (d,\delta)\in(0,1)$, such that
for $\gamma\in(0,\gamma_{0}]$, $h=D^{2}u$,  and  $\rho=R_{0}/\nu$,
 we have
\begin{equation}
                                   \label{12.15.3}
h^{\sharp}_{ \gamma,\rho}\leq
N\nu^{d/\gamma}\bM^{1/d}\big[|F[u]|^{d}\big]+N\tau_{0}\nu^{d/\gamma}
+N(\mu\nu^{d/\gamma}+\nu^{-\alpha})
\bM^{1/(\xi'd)}\big[|h|^{\xi'd}\big],
\end{equation}
where $\xi'=(\xi-1)/\xi$ and the constants $N$ depend
only on $d$, $K_{F}$,   and $\delta$.
\end{lemma}

We write $w\in A_{p}(\bR^{d})$ if $w$ is an $A_{p}$-weight on $\bR^{d}$.

\begin{lemma}
                          \label{lemma 12.15.4}
(i)
There exists $\gamma_{0}=\gamma_{0}(d,\delta)\in(0,1)$ such that for any $u\in W^{2}_{d,\loc}(\bR^{d})$, $\rho>0$, and $\gamma\in(0,\gamma_{0}]$ we have
\begin{equation}
                                   \label{12.15.6}
\bM^{1/\gamma}_{\rho}(|D^{2}u|^{\gamma})\leq N \bM^{1/d}_{\rho}(|F[u]|^{d})
+N   \rho^{-1}  \bM ^{1/d}_{ \rho}(|Du|^{d})
+N \rho^{-2} \bM^{1/d}_{\rho}(|u|^{d}),
\end{equation}
where the constants $N$ depend only on $d$,
 $\delta$,   and $K_F$.

(ii)
For any $\rho>0$, $p\in [1,\infty)$, and
$u\in W^2_{p,\loc}(\bR^{d})$, we have
\begin{equation}
                                            \label{12.15.7}
  \bM _{\rho}(|Du|^{ p })\leq N\bM^{1/2}_{\rho}
(|D^{2}u|^{ p })
\bM^{1/2}_{\rho}(| u|^{ p })+  N\rho^{-p}
\bM _{\rho}(| u|^{ p }),
\end{equation}
where the constants $N$ depend only on $d$ and $p$.

 (iii) For any $\rho>0$,  $K_{0}$, $p\in
(1,\infty)$, $w\in A_p(\bR^{d})$ with $[w]_{p}\le K_0$, and $u\in
W^2_{p,w}(\bR^{d})$, we have
\begin{equation}
                                            \label{1.38}
\int_{\bR^{d}}|Du|^{p}w\,dx\leq \rho^{p}
\int_{\bR^{d}}
\big |D^{2}u|^{p}w\,dx+N\rho^{-p}\int_{\bR^{d}}
\big |u|^{p}w\,dx,
\end{equation}
where $N$ depends only on $d$, $p$, and $K_0$.
\end{lemma}

Proof. First write $F[u]=a^{ij}D_{ij}u$ and take $r\geq\rho$
and a function $\zeta\in C^{\infty}_{0}(\bR^{d})$
such that $\zeta=1$ on $B_{r}$, $\zeta=0$ outside
$B_{2r}$, and
 $$
|D\zeta|\le N/r\le N/\rho,\quad |D^2\zeta|\le N/r^2\le N/\rho^2.
$$

Then
by  a result of  Fang-Hua Lin  \cite{Lin86},
$$
\dashint_{B_{r}}|D^{2}u|^{\gamma}\,dx
\leq \dashint_{B_{2r}}|D^{2}(\zeta  u)|^{\gamma}\,dx
$$
$$
\leq N\Big( \dashint_{B_{2r}}|\zeta  F[u]+a^{ij}2D_{i}\zeta
D_{j}u+ ua^{ij}D_{ij}\zeta )|^{d}\,dx\Big)^{\gamma/d}.
$$
This proves \eqref{12.15.6}.

The fact that  $\zeta=1$ on $B_r $   and multiplicative
inequalities show
that
$$
\dashint_{B_{r}}|D  u|^{ p }\,dx \leq N
\dashint_{B_{2r}}| D (\zeta u)|^{p}\,dx
$$
 $$
\leq N\Big(\dashint_{B_{2r}}| D^{2} (\zeta u)|^{p}\,dx\Big)^{1/2}
 \Big(\dashint_{B_{2r}}|u|^{p}\,dx\Big)^{1/2},
$$
where for $r\geq\rho$,
$$
\dashint_{B_{2r}}| D^{2} (\zeta u)|^{p}\,dx
\leq N\bM_{\rho}(|D^{2}u|^{p})+N  \rho^{-p} \bM _{\rho}(|D u|^{p})
+N  \rho^{-2p} \bM _{\rho} (|
u|^{p}).
$$
  Hence,
$$
 \bM _{\rho}(|D u|^{p})\leq N\Big(\bM_{\rho}(|D^{2}u|^{p})+
 \rho^{-p}\bM_{\rho}(|D
u|^{p})
\Big)^{1/2} \bM_{\rho}^{1/2}  (|  u|^{p})+N
\rho^{-p}\bM_{\rho} (|  u|^{p}),
$$ and \eqref{12.15.7} follows.

Finally, we prove \eqref{1.38}. We take an integer $m$ such
that $c2^{-m}\in (\rho/2,\rho]$. By  Remark
\ref{remark 5.27.1}, Corollary
\ref{corollary 12.3.1} with $\gamma=1$,  \eqref{12.15.5},
and the weighted Hardy-Littlewood
maximal function theorem (see more about this in the proof of Theorem
\ref{theorem 12.15.1})
$$
\int_{\bR^{d}}|Du|^{p}w\,dx\leq N
\int_{\bR^{d}}\big(  (Du)^{\#}_{1,m}
+\cM _{m}(|Du| )\big)^{p}w\,dx
$$
\begin{equation}
                        \label{eq1.50}
\leq N
\int_{\bR^{d}}\big(  (Du)^{\sharp}_{1,\rho}
+\bM _{\rho/2}(|Du| )\big)^{p}w\,dx.
\end{equation}
To apply Corollary
\ref{corollary 12.3.1} formally we need a certain
condition on the averages of $|Du|$. However, we always can use
cut-off functions and pass to the limit.
By Poincar\'e's inequality,
$$
(Du)^{\sharp}_{1,\rho}\leq N\rho \bM(|D^{2}u|).
$$
This together with \eqref{eq1.50} and \eqref{12.15.7} with $p=1$ gives
$$
\int_{\bR^{d}}|Du|^{p}w\,dx
\leq N
\int_{\bR^{d}}\big(  \rho\bM(|D^{2}u|)
+\rho^{-1}\bM _{\rho/2}(|u| )\big)^{p}w\,dx,
$$
which, by the weighted Hardy-Littlewood
 maximal function theorem, is bounded by the right-hand side
of \eqref{1.38}. The lemma is proved.

Estimate \eqref{1.38} admits the following
localization.
 \begin{lemma}
                                      \label{lemma 4.10.10}

For any $\rho\in (0,\infty)$, $\varepsilon\in (0,1]$,
$K_{0}$,
$p\in (1,\infty)$, $w\in A_p(\bR^{d})$ with $[w]_{p}\le K_0$, and
$u\in W^2_{p,w}(B_\rho)$, we have
\begin{equation}
                                            \label{2.11}
\int_{B_{\rho/2}}|Du|^{p}w\,dx\leq \varepsilon \rho^p
\int_{B_\rho}
\big |D^{2}u|^{p}w\,dx+N \varepsilon ^{-1}
\rho^{-p}\int_{B_\rho}
\big |u|^{p}w\,dx,
\end{equation} where $N$ depends only on $d$, $p$, and
$K_0$.
\end{lemma}

Proof.
By scaling and noting that
$[w(\rho\cdot)]_p=[w]_p$, we may assume that $\rho=1$. For
$k=1,2,\ldots$, we take $\rho_k=1-2^{-k}$, $B^k=B_{\rho_k}$, and
$\zeta_k\in C_0^\infty(B^{k+1})$ such that $\zeta_k=1$ on
$B^k$ and
$$ |D\zeta_k|\le N2^k,\quad |D^2\zeta_k|\le N2^{2k}.
$$ It follows from \eqref{1.38} that for any
$\varepsilon_0 \in (0,1]) $
$$
\int_{B^k}|Du|^{p}w\,dx
\le \int_{B^{k+1}}|D(\zeta^k u)|^{p}w\,dx
$$
$$
\le \varepsilon_0 2^{-kp} \int_{\bR^d}|D^2(\zeta^k
u)|^{p}w\,dx +N\varepsilon_0^{-1} 2^{kp}
\int_{\bR^d}|\zeta^k u|^{p}w\,dx
$$
$$
\le N\varepsilon_0 \int_{B^{k+1}}|Du|^{p}w\,dx +N
\int_{B^{k+1}}\big(\varepsilon_0 2^{-kp}|D^2u|^{p}
+\varepsilon_0^{-1} 2^{kp}|u|^{p}\big)w\,dx.
$$
Now to get \eqref{2.11}, it suffices to multiply both
sides by $(N\varepsilon_0)^k$, sum in $k=1,2,\ldots$, and
take a sufficiently small $\varepsilon_0$ according to
$\varepsilon$.

\begin{lemma}
                                         \label{lemma 4.12.1}
In Lemma \ref{lemma 12.15.4} (iii) the condition
$u\in
W^2_{p,w}(\bR^{d})$ can be replaced with
$u\in
W^2_{p,w,\loc}(\bR^{d})$.

\end{lemma}

Proof. We may assume that the right-hand side of
\eqref{1.38}  is finite. In this situation
plug $B_{1}(x_{0})$ and  $B_{2}(x_{0})$
in place of $B_{ \rho /2}$ and $B_{ \rho }$ into \eqref{2.11}
with $\varepsilon= 1/2 $ and integrate with respect to $x_{0}$
over $\bR^{d}$. Then we will see that $Du\in L_{p,w}(\bR^{d})$.
After that Lemma \ref{lemma 12.15.4} (iii) yields the result.
The lemma is proved.

In the future we will use the following.
\begin{lemma}
                                         \label{lemma 4.10.1}
Let $0<r<R<\infty$, $\varepsilon\in (0,1]$,
 $K_{0}$,
$p\in (1,\infty)$, $w\in A_p(\bR^{d})$ with $[w]_{p}\le K_0$, and
$u\in W^2_{p,w}(B_{R})$. Then
$$
\int_{B_{r }}|Du|^{p}w\,dx\leq \varepsilon (R-r)^{p}
\int_{B_R}
\big |D^{2}u|^{p}w\,dx+N(\varepsilon (R-r))^{-p}\int_{B_R}
\big |u|^{p}w\,dx,
$$
where $N$ depends only on $d$, $p$, and
$K_0$.
\end{lemma}

This lemma is a simple corollary of Lemma
\ref{lemma 4.10.10}. Indeed, set $\rho= R-r$
  in Lemma
\ref{lemma 4.10.10} and plug $B_{\rho/2}(x_{0})$
and $B_{\rho }(x_{0})$ into \eqref{2.11}
in place of $B_{\rho/2}$
and $B_{\rho}$, respectively, with $x_{0}\in B_{r}$.
Then it will only remain to integrate the resulting
inequality with respect to $x_{0}$ over $B_{r}$.

\begin{remark}
                          \label{remark 4.22.1}
Below we use a few times the fact that
if $w$ is an $A_{p/ d }(\bR^{d})$-weight for some $p\in (d,\infty)$, then by definition
$w^{-1}\in L_{d/(p-d),\loc}(\bR^{d})$. Hence by H\"older's
inequality, if $f\in L_{p,w,\loc}(\bR^d)$,
then $|f|^{d}\in L_{1,\loc}(\bR^d)$.
In particular, if  $u\in W^2_{p,w,\loc }(\bR^{d})$, then
$|u|^d, |Du|^d, |D^2 u|^d \in L_{1,\loc}(\bR^d)$,
which implies that $u\in W^2_{d,\loc}(\bR^d)$.
\end{remark}

\begin{theorem}
                                      \label{theorem 12.15.1}

Take $R\in(0,\infty)$,
 $K_{0}\in(1,\infty)$.
Let   $p>   d$
 and let $w\in A_{p/ d }(\bR^{d})$
with $[w]_{p/d}
\leq K_{0}$.
Suppose that $D^{2}u \in L_{p,w}(\bR^d)$   and $u$
vanishes outside
$B_{R}$.  Then there exists
  $\theta=\theta(d,\delta,K_{F},  p,   K_{0})
\in(0,1)$
such that if Assumption \ref{assump1} ($\theta$)
is satisfied, then
$$
\int_{\bR^{d}}|D^{2}u|^{p}w\,dx
\leq N\int_{\bR^{d}}|F[u]|^{p}w\,dx
$$
\begin{equation}
                                        \label{12.15.2}
+N \int_{\bR^{d}}|u|^{p}w\,dx
+N\tau_{0}^{p}\int_{\bR^{d}}I _{B_{R+R_{0} }}w\,dx,
\end{equation}
where  $N$  is a constant
depending only on $d$, $\delta$, $K_F$, $K_0$,
$p$, and $R_0$.
\end{theorem}

Proof.
It is well known that the
appropriately stated
 Hardy-Littlewood maximal function
theorem holds for $A_{p}$-weights.
Therefore, by Remark
\ref{remark 5.27.1}, Corollary
\ref{corollary 12.3.1}, and \eqref{12.15.5}
the left-hand side of
\eqref{12.15.2} is less than a constant times
$$
\Big(\int_{\bR^{d}}|D^{2}u|^{p}w\,dx\Big)^{(p-\gamma\beta)/p}
\Big(\int_{\bR^{d}}\big[(D^{2}u)^{\#}_{\gamma,m}
+\cM^{1/\gamma}_{m}\big(|D^{2}u|^{\gamma}\big)\big]^{p}w\,dx\Big)^{\gamma
\beta/p},
$$
where $\gamma\in (0,1)$ is a constant
depending only on $d$ and $\delta$
taken from Lemma \ref{lemma 12.15.1}. It follows that
the left-hand
side of \eqref{12.15.2} is less than a constant times
\begin{equation}
                                                    \label{12.15.4}
\int_{\bR^{d}}\big((D^{2}u)^{\#}_{\gamma,m}\big)^{p} w\,dx
+\int_{\bR^{d}}\cM^{p/\gamma}_{m}\big(|D^{2}u|^{\gamma}\big) w\,dx.
\end{equation}

By a reverse H\"older's inequality (also called
self-improving property  of $A_p$ weights, see for instance, Corollary 9.2.6 of \cite{Gr09}), we can find
$\xi\in (1,\infty)$ depending only on $d$, $p$, and $K_0$
such that $p>\xi' d$
($\xi'=\xi/(\xi-1))$,
$w\in A_{p/(\xi' d)}$, and $[w]_{p/(\xi' d)}\le N(d,K_0)$.
This is the first step to specify $\theta$
which will be taken from Lemma \ref{lemma 12.15.1}
after we find an appropriate $\mu>0$.
To this end, take $\nu\geq2$ to be specified
later and for $m$ such that $2^{-m}\sim R_{0}/\nu$
use \eqref{12.15.5} and
 \eqref{12.15.3}
 to estimate the first integral in \eqref{12.15.4}.
 Observe that $h^{\sharp}_{\gamma,\rho}$
vanishes outside $B_{R+R_{0}}$ and therefore we only need
to integrate the right-hand side of \eqref{12.15.3}
over this ball. This  gives  the last term
in \eqref{12.15.2} (after we fix $\nu$).

Then we again use the well-known properties
of $A_{p}$-weights mentioned above  and Lemma
\ref{lemma 12.15.1}
 to conclude that the first term in
\eqref{12.15.4}
is less than $\nu^{pd/\gamma}$ times the last term in
\eqref{12.15.2} plus
$$
N\nu^{pd/\gamma}\int_{\bR^{d}}|F[u]|^{p}w\,dx
+N(\mu\nu^{d/\gamma}+\nu^{-\alpha})^{p}
\int_{\bR^{d}}|D^{2}u|^{p}w\,dx.
$$
We choose first large $\nu$ and then small $\mu$
to absorb the last expression, which is {\em finite\/},
into the left-hand side of \eqref{12.15.2}.
This shows how to choose $\mu$
and now we take $\theta$ from Lemma \ref{lemma 12.15.1}.
After that it only remains to use Lemma \ref{lemma 12.15.4}
in order to estimate the second term in \eqref{12.15.4}
first taking care of adjusting $\gamma=\gamma(d,\delta)$ to fit both Lemmas \ref{lemma 12.15.1} and \ref{lemma 12.15.4}.
The theorem is proved.

\begin{remark}
Scalings show that the only constant $N$
in \eqref{12.15.2} depending on $R_{0}$
is the one in front of the integral of $|u|^{p}w$.
This one equals $N(d,\delta,K_{F},  p,   K_{0})R_{0}^{-2p}$.

\end{remark}

\begin{lemma}
                                    \label{lemma 4.8.1}
Let $u\in W^{2}_{d,\loc}(\bR^{d})$
be bounded and let $p>d$. Then
for any $\bS_{\delta}$-valued function
$a$ on $\bR^{d}$
$$
|u|^{p}\leq N(\delta,d,p)\bM(|a^{ij}D_{ij}u-u|^{p}).
$$
In particular, under Assumption \ref{assumption 4.26.1}
$$
|u|^{p}\leq N(\delta,d,p)\bM(|F[u]-u|^{p}).
$$
\end{lemma}

Proof.  First observe that the second estimate
follows from the first one since
 $ F[u] =a^{ij}D_{ij}u $, where
$(a^{ij})$ is an appropriate $\bS_{\delta}$-valued
function.
To prove the first estimate,
let $G(x,y)$ be a Green's function of
$L:= a^{ij}D_{ij}-1$ in $\bR^{d+1}$ and $f=-Lu$. Then
  we have
$$
u(0)= \int_{\bR^{d}}
G( 0,y)f( y)\,dy .
$$
Hence,
$$
|u(0)|\leq \int_{\bR^{d}}
G( 0,y)|f( y)|\,dy.
$$

We are going to use the following estimate
easily obtained, say, by probabilistic arguments:
for any $\beta\geq0$
\begin{equation}
                                          \label{4.9.1}
\int_{\bR^{d } }G(0,y) |y|^{\beta}
\,dy \leq N(\alpha,d,\delta).
\end{equation}

Observe that for any $h( y)\geq0$ and $\alpha>0$
$$
\int_{1}^{\infty}r^{-\alpha-1}
\Big(\int_{B_{r}}h(y) (
|y|^{\alpha}\vee1)\,dy \Big)dr
$$
$$
=\int_{\bR^{d } }h(y)(
|y|^{\alpha}\vee 1)\Big(\int_{ |y|\vee 1 }^{\infty}
r^{-\alpha-1}\,dr\Big)
dy
=\frac{1}{\alpha}\int_{\bR^{d } }h \,
dy .
$$
Furthermore, by using \eqref{4.9.1}, H\"older's inequality,
and the   Aleksandrov estimate, for  $q=p/d>1$,
we get
$$
\int_{B_{r}}G(0,y)|f(y)|(
|y|^{\alpha}\vee 1)\,dy \leq N
\Big(\int_{B_{r}}G(0,y)|f(y)|^{q}
\,dy \Big)^{1/q}
$$
$$
\leq N\Big(\int_{B_{r}}|f( y)|^{p}
\,dy \Big)^{1/p}\leq Nr^{d/p}\big(\bM(|f|^{p})(0)\big)^{1/p}.
$$
For $d/p-\alpha-1<-1$ we get the desired result by
integrating in $r\in [1,\infty)$ and collecting the above
estimates.  The lemma is proved.

Thanks to the properties
of $A_{p}$-weights mentioned in the beginning of the
proof of Theorem \ref{theorem 12.15.1},
  we have the following.

\begin{corollary}
                                       \label{corollary 4.25.1}
Under the assumption of Lemma \ref{lemma 4.8.1}
take $K_{0}\in(1,\infty)$ and let $w$ be an $A_{p/d}$-weight
with $[w]_{p/d}\leq K_0$. Then there exists a constant
$N=N(\delta,d,p,K_{0})$ such that
\begin{equation}
                                        \label{4.24.6}
\int_{\bR^{d}} |u|^{p} w\,dx
\leq N\int_{\bR^{d}}|a^{ij}D_{ij}u-u|^{p}w\,dx,
\end{equation}
in particular, under Assumption \ref{assumption 4.26.1}
$$
\int_{\bR^{d}} |u|^{p} w\,dx
\leq N\int_{\bR^{d}}|F[u]-u|^{p}w\,dx.
$$
\end{corollary}

Here is a generalization for fully nonlinear
operators
with almost VMO dependence on $x$
considered  in Sobolev spaces with $A_{p}$-weights
of the classical Sobolev estimates known for linear
operators with continuous coefficients.

\begin{theorem}
                                      \label{theorem 4.8.1}
Take $K_{0}\in(1,\infty)$.
Let $p>d$ and let $w\in A_{p/ d }(\bR^{d})$
with $[w]_{p/d}
\leq K_{0}$.
Let
\begin{equation}
                                          \label{4.11.1}
u\in   W^2_{p,w }(\bR^{d}).
\end{equation}

  Then there exists
  $\theta=\theta(d,\delta, K_{F},   p,   K_{0})
\in(0,1)$
such that, if Assumption \ref{assump1} ($\theta$)
is satisfied with $\tau_{0}=0$, then
\begin{equation}
                                               \label{4.8.1}
\int_{\bR^{d}}(|D^{2}u|^{p}+|Du|^{p})w\,dx
\leq N\int_{\bR^{d}}|F[u]|^{p}w\,dx
+N \int_{\bR^{d}}|u|^{p}w\,dx,
\end{equation}
and if, in addition, $u$ is bounded and
Assumption \ref{assumption 4.26.1} is satisfied
then
\begin{equation}
                                               \label{4.8.2}
\int_{\bR^{d}}(|D^{2}u|^{p}+|Du|^{p}+|u|^{p})w\,dx
\leq N\int_{\bR^{d}}|F[u]-u|^{p}w\,dx,
\end{equation}
where  the constants $N$
depend  only on $d$, $\delta$, $K_F$, $K_0$,
$p$, and $R_0$.
\end{theorem}

Proof. To prove \eqref{4.8.1}, thanks to Lemma
 \ref{lemma 12.15.4} (iii), it suffices to estimate $|D^{2}u|$.

To this end introduce
$\zeta\in C^{\infty}_{0}(\bR^{d})$ such that $\zeta(0)=1$
and plug $u_{n} :=u \zeta_{n}$, where $\zeta_{n}(x)=
\zeta(x/n)$, into \eqref{12.15.2} (just in case remember that \eqref{12.15.2}
is proved for functions with compact support). Then the result
 follows by the
dominated convergence theorem from the fact that
$$
|F[u]-F[u_{n}]|\leq K_{F}|D^{2}u-D^{2}u_{n}|
$$
$$
\leq
N|1-\zeta_{n}|\,|D^{2}u|+Nn^{-1}|Du|+Nn^{-  2 }|u|.
$$

To prove \eqref{4.8.2} it suffices to  use Corollary \ref{corollary 4.25.1}.
The theorem is proved.

\begin{theorem}
                            \label{theorem 4.9.2}
Let    $p_{i}>d$,  $i= 1, 2,\ldots,d$.
Assume that $u\in W^{2}_{1,\loc}(\bR^{d})$
and Assumption \ref{assumption 4.26.1} is satisfied.
Then there exists $\theta=\theta(d,\delta, d
, p_1 ,\ldots,p_d )\in(0,1)$ such that, if Assumption \ref{assump1} ($\theta$)
is satisfied  with $\tau_{0}=0$, then
\begin{equation}
                                      \label{4.9.5}
\|   D^{2}u,Du,u \|_{L_{p_{1},\ldots,p_{d}}(\bR^{d})}
\leq N\| F[u]-u \,\|_{L_{p_{1},\ldots,p_{d}}(\bR^{d})}
\end{equation}
provided  that the left-hand side is finite,
where
$$
\| f
\|_{L_{p_{1},\ldots,p_{d}}(\bR^{d})}^{p_{d}}
$$
\begin{equation}
                                       \label{4.9.6}
:=\int_{\bR  }\Big(\cdots
\Big(\int_{\bR }
\Big(\int_{\bR } |f|^{p_{1}}\,dx_{1}
\Big)^{p_{2} /p_1 }\,dx_{2}\Big)^{p_{3}/p_{2}}
\cdots\Big)^{p_{d} /p_{d-1} }\,dx_{d},
\end{equation}
and the constant  $N$ depends only on $d$, $\delta$, $p_1,\ldots,p_d$, and $R_0$.
\end{theorem}

Proof. First we assume that $u$ is smooth and
has compact support.  Then \eqref{4.9.5} follows from
 Theorems \ref{theorem 4.8.1} and \ref{theorem 12.2.1}.

If $u$ is just smooth (and the left-hand side
of \eqref{4.9.5} is finite), one can use the
same approximation of $u$
as in the proof of \eqref{4.8.1}.

Finally, in the general case introduce
$u^{(\varepsilon)}$ as the mollified $u $.
By the Minkowski inequality (the norm of a sum is less than
the sum of norms)
the left-hand side
of \eqref{4.9.5} with $u^{(\varepsilon)}$
in place of $u$
 is less than its original.
After that writing \eqref{4.9.5}
with $u^{(\varepsilon)}$
in place of $u$, using the Lipschitz continuity
of $F(\sfu'',x)$ with respect to $\sfu''$,  noting $D^{2}u^{(\varepsilon)}
\to D^{2}u$ in the above mixed norm (see,
for instance, \cite{BP61}),
and letting $\varepsilon\downarrow 0$, we easily
finish the proof.
The theorem is proved.

Sometimes in the sequel we consider
$F$'s that are positive homogeneous in $\sfu''$.
In that case we impose the following.

 \begin{assumption}[$\theta$]
                                \label{assump1homo}
  (i) The function $F$ is Lipschitz continuous with respect to $\sfu''$
with Lipschitz constant $K_{F}$ and
is positive homogeneous of degree one with respect to $\sfu''$.

There exists $R_0\in(0,1]$ such that, if   $r\in (0, R_0]$ and $z\in \bR^{d}$, then one can find a {\em convex\/} function $\bar{F} (\sfu'' )=
\bar{F}_{z,r } (\sfu'' )$ (independent
of $x $) positive homogeneous of degree one,  for which

(ii) We have $ D_{\sfu'' }\bar{F} \in\bS_{\delta}$
at all points of differentiability
of $\bar{F}$;

(iii)
For any $\sfu''\in\bS$ with $|\sfu''|=1$,
 we have
\begin{equation}
                                                \label{7.30.2homo}
\int_{ B_{r}(z)}\big|F\big(\sfu'' ,x\big)-\bar{F}( \sfu'')\big| \,dx\leq
\theta
\big| B_{r}(z)\big|.
\end{equation}
\end{assumption}

\begin{remark} It is worth noting that if $F$ is positive
homogeneous of degree one with respect to $\sfu''$ and
satisfies Assumption \ref{assump1}, then it also satisfies
Assumption \ref{assump1homo}. Indeed, let $\bar F$ be the
function from Assumption \ref{assump1}. Then the function
$\limsup_{\lambda\to \infty} \lambda^{-1} \bar F(\lambda
u'')$ is convex, positive homogeneous of degree one, and
satisfies Assumption \ref{assump1homo} (ii) and (iii).
\end{remark}

Sometimes the following result is useful.
\begin{lemma}
                                      \label{lemma 12.16.1b}

Take $R\in(0,\infty)$, $K_{0}\in(1,\infty)$, $p>d$ and let
$w\in A_{p/ d}(\bR^{d})$  with $[w]_{p/d}\leq K_{0}$. Suppose that
\begin{equation}
                                          \label{2.20.1}
u\in W^{2}_{p,w}(B_R).
\end{equation}
  Then there exists
  $\theta=\theta(d,\delta,K_{F}, p,   K_{0})
\in(0,1) $
such that, if Assumption  \ref{assump1homo}  ($\theta$)
is satisfied, then for any $r\in(0,R)$
$$
\int_{B_{r}}|D^{2}u|^{p}w\,dx
\leq N\int_{B_{R}}|F[u]|^{p}w\,dx
$$
\begin{equation}
                                                    \label{12.15.2b}
+N \int_{ B_R }((R-r)^{-1}|Du|+((R-r)^{-2} +1 )|u|)^{p}w\,dx,
\end{equation}
where  $N$  is a constant depending
only on $d$, $\delta$, $K_F$, $K_0$, $p$, and $R_0$.
\end{lemma}

Proof. Take  a nonnegative  $\zeta\in C^{\infty}_{0}(B_{R})$ such that
$\zeta=1$ on $B_{r}$ and $|D\zeta|\leq N(R-r)^{-1}$,
$|D^{2}\zeta|\leq N(R-r)^{-2}$.
 It follows from \eqref{2.20.1}
that $\zeta u\in W^2_{p,w}(B_R)$.
Then apply Theorem \ref{theorem 12.15.1}
to $\zeta u$ after observing that due to the homogeneity
of $F$ we have
$$
|\zeta F[u]-F[\zeta u]|\leq N(|D\zeta|\,|Du|+|D^{2}\zeta|\,|u|).
$$
This yields the result.

The following result and Lemma \ref{lemma 4.8.1}
easily imply Theorem \ref{theorem 4.8.1} once more,
however in Theorem \ref{theorem 4.8.1} we do not assume
that $F$ is positive homogeneous.
\begin{lemma}
                                      \label{lemma 4.10.6}
Let the assumptions of Lemma
\ref{lemma 12.16.1b}
be satisfied and take $\theta$ from that lemma.
 Then there is a constant $N$ depending
only on $d$, $\delta$, $K_F$, $K_0$, $p$, and $R_0$,
such that for any $r\in  [R-1 ,R)$, $r>0$,
\begin{equation}
                                             \label{4.10.1}
\int_{B_{r}}|D^{2}u|^{p}w\,dx
\leq NF
+N(R-r)^{-2p}U,
\end{equation}
\begin{equation}
                                             \label{4.11.3}
\int_{B_{r}}|D u|^{p}w\,dx
\leq N(R-r)^{p}F
+N(R-r)^{-p}U,
\end{equation}
where
$$
F=\int_{B_{R}}
|F[u]|^{p}w\,dx,\quad U=\int_{B_{R}}
|u|^{p}w\,dx.
$$
\end{lemma}

Proof. For $k=0,1,\ldots$ set $\rho_{k}=R-2^{-k}(R-r)$,
$B^{k}=B_{\rho_{k}}$ and find $\zeta_{k}\in C^{\infty}_{0}
(B^{k+1})$ such that $\zeta_{k}=1$ on $B^{k}$
and $|D\zeta_{k}|\leq N(R-r)^{-1}2^{k}$,
$|D^{2}\zeta_{k}|\leq N(R-r)^{-2}2^{2k}$, where
$N=N(d)$.

By Lemma \ref{lemma 12.16.1b} we have
$$
D''_{k}:=\int_{B^{k}}|D^{2}u|^{p}w\,dx \leq
NF+N(R-r)^{-p}2^{kp}D'_{k+1}+N(R-r)^{-2p}2^{2kp}U,
$$
where
$$
D'_{k}=\int_{B^{k}}|Du|^{p}w\,dx.
$$
By Lemma \ref{lemma 4.10.1} for $\varepsilon\in(0,1]$
\begin{equation}
                                            \label{4.11.2}
D'_{k}\leq \varepsilon 2^{-kp}(R-r)^{p}D''_{k+1}
+N\varepsilon^{-1}(R-r)^{-p}2^{kp}U.
\end{equation}
It follows that
$$
D''_{k}\leq NF+\varepsilon N_{1}D''_{k+2}+
N\varepsilon^{-1}(R-r)^{-2p}2^{2kp}U
$$
We choose $\varepsilon$ so that $\varepsilon N_{1}\leq 2^{-6p}$,
multiply this inequality by $2^{- 3kp}$ and sum up with respect
to even $k$ from $0$ to $\infty$. Then we cancel like terms
(which are finite since $D^{2}u\in L_{p,w}(B_{R})$)
and come to \eqref{4.10.1}.

After that \eqref{4.11.3} follows from \eqref{4.10.1}
and \eqref{4.11.2}. The lemma is proved.

By substituting $B_{r}(x_{0})$
and $B_{R}(x_{0})$ in place of $B_{r} $
and $B_{R} $, respectively, then
taking $R=2r=1$ and integrating with respect to $x_{0}$
over $\bR^{d}$ we obtain the following.
\begin{theorem}
                                      \label{theorem 4.11.1}
Theorem \ref{theorem 4.8.1} remains true
if condition \eqref{4.11.1} is replaced with
$u\in   W^2_{p,w,\loc}(\bR^{d})$ provided,
additionally,
that $F$ is positive homogeneous
of degree one with respect to $\sfu''$.
In particular, if $u$ is bounded and
 the right-hand side of
\eqref{4.8.2} is finite, then $u\in W^{2}_{p,w}(\bR^{d})$.

\end{theorem}
\begin{remark}
                                      \label{remark 4.11.1}

Generally,  \eqref{4.8.2} may fail if $u$ is unbounded.
Indeed, if $d=1$ and $F[u]=u''$, the function $e^{x}$
satisfies $F[u]-u=0$ and is nonzero.
\end{remark}

\begin{remark}
                                    \label{remark 2.20.1}
Condition \eqref{2.20.1} is not well suited
for application of
the extrapolation theorem
of J. L. Rubio de Francia \cite{MR745140}.
 In this connection it is useful to know
that, for any $K_{0},p\in(1,\infty)$ and $w\in A_{p}
  (\bR^{d})$
with $[w]_{p}\leq K_{0}$ there exists
$q=q(d,K_{0},p)\in(1,\infty)$ such that
$W^{2}_{q}(B_{R})\subset W^{2}_{p,w}(B_{R})$
for any $R<\infty$.
This follows from the fact
that (see, for instance, Corollary 9.2.4 of \cite{Gr09})  $w$
is in $L_{r,\loc}(\bR^{d})$ for
an appropriate $r >1 $
 depending only on $d$, $p$, and $K_0$.
\end{remark}

\begin{theorem}
                                   \label{theorem 12.29.1}

Take $R\in(0,\infty)$, $r\in (0,R)$,   $p>  d$  and
take  $p_{i}>d$,  $i= 1, 2,\ldots,d$.
Assume that $u\in W^{2}_{p}(B_{R})$.
Then there exists
$$
\theta=\theta(d,\delta, p ,p_{1},\ldots,p_{d})
\in(0,1)
$$
such that, if Assumptions  \ref{assump1homo}  ($\theta$)
 and \ref{assumption 4.26.1} are satisfied, then
$$
\| I_{B_{r}}  D^{2}u,  I_{B_{r}} Du \|_{L_{p_{1},\ldots,p_{d}}(\bR^{d})}
\leq N\|I_{B_{R} } F[u] \,\|_{L_{p_{1},\ldots,p_{d}}(\bR^{d})}
$$
\begin{equation}
                                      \label{2.20.3}
+N\|I_{B_{R} } u\|_{L_{p_{1},\ldots,p_{d}}(\bR^{d})},
\end{equation}
where the constants
 $N$ depend only on  $r$, $R$,
$d$, $\delta$, $p_1,\ldots,p_d$, and $R_0$.
\end{theorem}

Proof. In Theorem
\ref{theorem 12.2.1} take $m=d$, $K_{0}=1$, $k(1)=\ldots.=k(d)=1$
and take $\Lambda_{0}$ from there which now
depends only on $d$ and $p_{1},\ldots,p_{d}$.
Then take $q=q(d,\Lambda_{0},p_{1})$ from
Remark \ref{remark 2.20.1} and assume that
$u\in W^{2}_{q}(B_{R})$. In that case in light of
Remark \ref{remark 2.20.1} estimate
\eqref{2.20.3} follows from Lemma
\ref{lemma 4.10.6} and Theorems \ref{theorem 12.2.1}.

In the general case, we may assume that the right-hand side
of \eqref{2.20.3} is finite and introduce $f=F[u]$
and $f^{(\varepsilon)}$ as the mollified $fI_{B_{R}}$.
By Minkowski's inequality (the norm of a sum is less than
the sum of norms) the above mixed norm of $f^{(\varepsilon)}$
is less than that of $fI_{B_{R}}$. Then
for small $\varepsilon >0$
 define smooth $u_{\varepsilon}$ so that they converge to $u$
uniformly on $\partial B_{R}$ and define
$u^{\varepsilon}$ as unique $W^{2}_{p
 }(B_{R }) $-solutions of $F[u^{\varepsilon}]
=f^{(\varepsilon)}$ in $B_{R }$ with boundary condition
$u^{\varepsilon}=u_{\varepsilon}$ on
$\partial B_{R }$.
 Such solutions exist
and belongs to $W^{2}_{  q }(B_{R })$
  thanks to Theorem
2.1 of \cite{Kr_13} (provided that
an appropriate choice of $\theta$ is made).

Owing to the Aleksandrov estimate,  $u^{\varepsilon}\to u$
uniformly on $B_{R }$ as $\varepsilon\downarrow 0$.
In light of \eqref{2.20.3}
 the mixed norms of
$D^{2}u^{\varepsilon}$ and $D u^{\varepsilon}$
are bounded and since $u^{\varepsilon}\to u$,
they weakly converge in the space with mixed norm to
$D^2u$ and $D u $. The norm of the weak limit is less than
the limit of norms and this proves
\eqref{2.20.3} and the theorem.

\mysection{Elliptic equations in half spaces. First approach}

                                    \label{section 4.10.1}

Here we consider elliptic equations in the
half-space $$
\bR^d_+:=\{x=(x_{1},x'):x_1\geq 0,
x'\in\bR^{d-1}\}
$$
{\em without\/}
boundary conditions, and prove estimates near the boundary
with $A_{p}$-weights on $\bR^d_+$. A
typical and probably the most
interesting example of $A_{p}$-weights   on $\bR^d_+$
is the distance to the boundary to some
power, i.e., $w(x)=x_1^q$. It is easy to see that %%%
$w\in
A_p (\bR^d_+)$ (that is, $w$ is an $A_{p}$-weight
on $\bR^{d}_{+}$)
 if and only if $q\in (-1,p-1)$.
The way to build our estimates is taken from
\cite{Kr_08_1}.

Our underlying $\Omega$ is $\bR^{d}_{+}$ and $\bC_{n}$
are the cubes from Section \ref{section 4.10.3}
only lying in $\bR^{d}_{+}$.
Naturally, $\bL$ is the set of continuous functions
on $\bR^{d}_{+}$ with compact support.

For $n\in\bZ$, $R>0$
introduce
$$
S_{n}=[2^{-n},2^{-n+1}]\times\bR^{d-1}
,\quad T_{n}=[2^{-n-1},2^{-n+2}]\times\bR^{d-1},
\quad
B^{+}_{R}=B_{R}\cap \bR^{d}_{+}.
$$

In this section we consider
   a function $F(\sfu'',x)$, $\sfu''\in\bS$,
$x\in  \bR^{d} $, that is positive homogeneous
of degree one with respect to $\sfu''$.

\begin{lemma}
                                    \label{lemma 4.8.2}

Take
 $K_{0}\in(1,\infty)$,
    $p>   d$,
 and let  $w\in A_{p/ d }(\bR^{d}_{+})$ with $[w]_{p/d}
\leq K_{0}$.
Let $u$ be a bounded function on $\bR^{d}_{+}$ such that
$$
u\in   W^2_{p,w}(\bR^{d}_{+}).
$$

  Then there exists
  $\theta=\theta(d,\delta,   p,   K_{0})
\in(0,1)$
such that if Assumption  \ref{assump1homo} ($\theta$)
is satisfied, then there is a constant $N$, depending only on $d$, $\delta$,  $K_0$,
$p$, and $R_0$, such that for any $n\in\bZ$
and any $\varepsilon \in (0,1]$ we have
$$
\int_{S_{n}}|D^{2}u|^{p}w\,dx
\leq N\int_{T_{n}}|F[u]-u|^{p}w\,dx
$$
\begin{equation}
                                          \label{4.8.4}
+N2^{pn}\int_{T_{n}}|D u|^{p}w\,dx
+N(2^{2pn}+1)\int_{T_{n}}| u|^{p}w\,dx,
\end{equation}
$$
\int_{S_{n}}|Du|^{p}w\,dx
\leq N\varepsilon 2^{-pn}\int_{T_{n}}|D^{2}u|^{p}w\,dx
$$
\begin{equation}
                                          \label{4.8.7}
+N\varepsilon\int_{T_{n}}|Du|^{p}w\,dx
+N\varepsilon^{-1}2^{pn}\int_{T_{n}}| u|^{p}w\,dx.
\end{equation}

Furthermore, for any $\varepsilon\in (0,1]$
$$
\int_{\bR^{d},x_{1}\geq 2}|D^{2}u|^{p}w\,dx
\leq N\int_{\bR^{d},x_{1}\geq 1}|F[u]-u|^{p}w\,dx
$$
\begin{equation}
                                          \label{4.8.01}
+N \int_{\bR^{d},x_{1}\geq 1}|D u|^{p}w\,dx
+N \int_{\bR^{d},x_{1}\geq 1}| u|^{p}w\,dx,
\end{equation}
$$
\int_{\bR^{d},x_{1}\geq 2}|Du|^{p}w\,dx
\leq N\varepsilon  \int_{\bR^{d},x_{1}\geq 1}|D^{2}u|^{p}w\,dx
$$
\begin{equation}
                                          \label{4.8.02}
+N\varepsilon\int_{\bR^{d},x_{1}\geq 1}|Du|^{p}w\,dx
+N\varepsilon^{-1} \int_{\bR^{d},x_{1}\geq 1}| u|^{p}w\,dx.
\end{equation}
\end{lemma}

Proof. To prove \eqref{4.8.4} we use the fact that
there is a nonnegative
 $\zeta\in C^{\infty}_{0}(\bR)$ such that
$\zeta=1$ on $[2^{-n},2^{-n+1}]$, $\zeta=0$
outside $[2^{-n-1},2^{-n+2}]$ and $2^{-n}|\zeta'|,
2^{-2n}|\zeta''|\leq N$, where $N$ is an absolute constant.
Then we
apply Theorem \ref{theorem 4.8.1}
to $u\zeta$ and, after observing that,
due to the positive homogeneity  and the Lipschitz continuity
of $F$, we have
$$
|\zeta F[u]-F[\zeta u]|
\leq N(|D\zeta|\,|Du|+|D^{2}\zeta|\,|u|),
$$
immediately arrive at \eqref{4.8.4}.
Of course, since we used a result in which
$w$ is an $A_{p/d}$-weight on $\bR^{d}$ rather than
on $\bR^{d}_{+}$, we first extend $w$
in an even way across $\{x_{1}=0\}$ with its norm controlled by $K_0$.
To prove \eqref{4.8.7} we use the same substitution
but into \eqref{1.38} and choose $\rho^{p}=
\varepsilon 2^{-pn}$.
Similarly \eqref{4.8.01} and \eqref{4.8.02} are obtained.
The lemma is proved.

\begin{theorem}
                                 \label{theorem 4.8.3}
Let $q\in \bR$.
Under the assumptions of Lemma \ref{lemma 4.8.2}
and for $\theta$ from that lemma,
if Assumption \ref{assump1} ($\theta$)
is satisfied,
then
$$
\int_{\bR^{d}_{+}}\hat x_{1}^{q}|\hat x_{1} D^{2}u|^{p}
w\,dx+\int_{\bR^{d}_{+}}\hat x_{1}^{q}|  D u|^{p}
w\,dx
$$
\begin{equation}
                                       \label{4.8.5}
\leq N
\int_{\bR^{d}_{+}}\hat x_{1}^{q}|\hat x_{1} (F[u]-u) |^{p}
w\,dx
+N\int_{\bR^{d}_{+}}\hat x_{1}^{q}|\hat x_{1}^{-1}  u|^{p}
w\,dx,
\end{equation}
 where $\hat x_1=\min\{x_1,1\}$,
provided that the left-hand side is finite,
where the $N$'s depend only on
 $d$, $\delta$,  $K_0$, $p$, $q$, and $R_0$.
\end{theorem}

Proof. Multiply both parts of \eqref{4.8.4}
 by $2^{-qn-pn}$,
sum up over $n\geq0$, and use the fact that
$2^{-qn-pn}\sim x_{1}^{q+p}$ on $S_{n}$ and $T_{n}$.
Then we get
$$
\int_{\bR^{d}_{+},x_{1}\leq2}x_{1}^{q}|x_{1}  D^{2} u|^{p}
w\,dx\leq N\int_{\bR^{d}_{+},x_{1}\leq 4}
x_{1}^{q}|x_{1}( F[u]-u)|^{p}
w\,dx
$$
$$
+N\int_{\bR^{d}_{+},x_{1}\leq 4}
x_{1}^{q}|   D  u|^{p}
w\,dx+N\int_{\bR^{d}_{+},x_{1}\leq 4}
x_{1}^{q}|x_{1}^{-1} u|^{p}
w\,dx.
$$

Multiplying \eqref{4.8.7} by $2^{-qn}$
and summing up yields
$$
\int_{\bR^{d}_{+},x_{1}\leq2}x_{1}^{q}|  D  u|^{p}
w\,dx\leq N\varepsilon\int_{\bR^{d}_{+},x_{1}\leq 4}
x_{1}^{q}|x_{1}D^{2}u|^{p}
w\,dx
$$
$$
+N\varepsilon\int_{\bR^{d}_{+},x_{1}\leq 4}
x_{1}^{q}|   D  u|^{p}
w\,dx+N\varepsilon^{-1}\int_{\bR^{d}_{+},x_{1}\leq 4}
x_{1}^{q}|x_{1}^{-1} u|^{p}
w\,dx.
$$

By combining  these estimates
 with \eqref{4.8.01} and \eqref{4.8.02}
we see that for any $\varepsilon\in(0,1]$
$$
\int_{\bR^{d}_{+} }\hat x_{1}^{q}|\hat x_{1}  D^{2} u|^{p}
w\,dx\leq N\int_{\bR^{d}_{+} }
\hat x_{1}^{q}|\hat x_{1}( F[u]-u)|^{p}
w\,dx
$$
\begin{equation}
                                       \label{4.9.7}
+N\int_{\bR^{d}_{+} }
\hat x_{1}^{q}|   D  u|^{p}
w\,dx+N\int_{\bR^{d}_{+} }
\hat x_{1}^{q}|\hat x_{1}^{-1} u|^{p}
w\,dx,
\end{equation}
$$
\int_{\bR^{d}_{+} }\hat x_{1}^{q}|  D  u|^{p}
w\,dx\leq N\varepsilon\int_{\bR^{d}_{+} }
\hat x_{1}^{q}|\hat x_{1}D^{2}u|^{p}
w\,dx
$$
\begin{equation}
                                       \label{4.9.8}
+N\varepsilon\int_{\bR^{d}_{+} }
\hat x_{1}^{q}|   D  u|^{p}
w\,dx+N\varepsilon^{-1}\int_{\bR^{d}_{+} }
\hat x_{1}^{q}|\hat x_{1}^{-1} u|^{p}
w\,dx.
\end{equation}
By choosing $\varepsilon$ in an obvious way,
we arrive at \eqref{4.8.5}. The theorem is proved.

 The next theorem follows from Theorems \ref{theorem 4.8.3}
and \ref{theorem 12.2.1}.

\begin{theorem}
                                  \label{theorem 4.10.2}
Let $p_{1},p_{2}>d$, $q\in \bR$, and let $u\in C^{\infty}_{0}(\bR^{d}_{+})$
have (closed) support in $\{x_{1}>0\}$. Then
there exists
  $\theta=\theta(d,\delta, q, p_{1},p_{2} )
\in(0,1)$
such that if Assumption  \ref{assump1homo} ($\theta$)
is satisfied, then there is a constant
$N$, depending only on $d$, $\delta$,
$q$, $p_{1}$, $p_{2}$,
  and $R_0$, such that
$$
\int_{0}^{\infty}\hat x_{1}^{q }
\Big(\int_{\bR^{d-1}}\Big[|\hat x_{1} D^{2}u|+|D^{2}u|
\big]^{p_{1}}\,dx'
\Big)^{p_{2}/p_{1}}\,dx_{1}
$$
$$
\leq N\int_{0}^{\infty}\hat x_{1}^{q}
\Big(\int_{\bR^{d-1}}|\hat x_{1}(F[u]-u)|^{p_{1}}\,dx'
\Big)^{p_{2}/p_{1}}\,dx_{1}
$$
\begin{equation}
                                              \label{4.11.4}
+N\int_{0}^{\infty}\hat x_{1}^{q}
\Big(\int_{\bR^{d-1}}|\hat x_{1}^{-1}u|^{p_{1}}\,dx'
\Big)^{p_{2}/p_{1}}\,dx_{1}.
\end{equation}
\end{theorem}

The reader understands that similar estimate holds
for mixed norms when we integrate with respect to
$x_{1}$ first.

\begin{remark}
                                   \label{remark 4.11.3}
Introduce a Banach space of functions on $\bR^{d}_{+}$
having finite norm defined by
$$
\|u\|^{p_{2}}=
\int_{0}^{\infty}\hat x_{1}^{q }
\Big(\int_{\bR^{d-1}}\Big[|\hat x_{1} D^{2}u|+|Du|+
\hat x_{1}^{-1}|u|\big]^{p_{1}}\,dx'
\Big)^{p_{2}/p_{1}}\,dx_{1}.
$$
It turns out that the set of
$u\in C^{\infty}_{0}(\bR^{d}_{+})$ that
have (closed) support lying in $\{x_{1}>0\}$
is everywhere dense in this space, so that
estimate \eqref{4.11.4} automatically extends
to all functions in this space.

To prove this, first take a smooth function
$\eta(r)$ such that $\eta(r)=0$ for $r<-1$ and $\eta(r)=1$
for $r>0$ introduce $\eta_{k}(x)=\eta(k^{-1}\ln x_{1})$,
$u_{k}=u\eta_{k}$ and by using the dominated convergence
theorem prove that, if $\|u\|<\infty$, then $\|u-u_{k}\|
\to 0$ as $k\to \infty$. After that it only remains
to apply usual tools to approximate $u_k$ by smooth functions
 which have (closed) support lying in $\{x_{1}>0\}$.

\end{remark}

In  the next section we show that for some values of $q$
it is possible to eliminate the last term in
\eqref{4.11.4}.

\mysection{elliptic equations in half spaces.
Second approach}
                                    \label{section 4.10.2}

We use the setting and the notation
from the beginning of Section
\ref{section 4.10.1} and in this section we deal
with a function $F(\sfu'',x)$ given for $\sfu''\in\bS$ and
 $x\in \bR^d_+$ and satisfying one of the following assumptions before which we
introduce
$$
B_r^+(x)=B_r(x)\cap \bR^d_+.
$$

\begin{assumption}[$\theta$]
                                       \label{assump1b}
Assumption \ref{assump1} ($\theta$) is satisfied if we
replace there  $\bR^{d}$ and $B_{r}(z)$    with $\bR^{d}_{+}$ and
$B_{r}^+(z)$, respectively.
\end{assumption}
 \begin{assumption}[$\theta$]
                                       \label{assump1bhomo}
Assumption \ref{assump1homo} ($\theta$) is satisfied if we
replace there  $\bR^{d}$ and $B_{r}(z)$  with $\bR^{d}_{+}$ and
$B_{r}^+(z)$, respectively.
\end{assumption}

Similarly, we introduce
$$
h^{\sharp}_{ \gamma,\rho}(x),\quad \bM_\rho h(x),\quad
 \text{and}\quad\bM h(x)
$$
on $\bR^{d}_{+}$ (by taking $B_{r}^+(x_{0})
\subset \bR^{d+1}_{+}$, $x_{0}\in \bR^{d}_{+}$).

From Lemma 4.2 of \cite{Kr_13} and the proof of Lemma 5.2 of
\cite{Kr_13}, we can easily obtain a boundary analog of Lemma
\ref{lemma 12.15.1}. This together with a boundary analog of
Lemma \ref{lemma 12.15.4} allows us to apply
Corollary \ref{corollary 12.3.1} and
 yields the following boundary estimate
corresponding to Theorem \ref{theorem 12.15.1} above.
\begin{theorem}
                                      \label{theorem 12.15.1b}

Take $R\in(0,\infty)$ and  $K_{0}\in(1,\infty)$. Let $p>d$ and let $w$ be an $A_{p/ d }$-weight on $\bR^d_+$ with $[w]_{p/d}
\leq K_{0}$.
Suppose that $D^{2}u \in L_{p,w}(\bR_+^d)$ and $u$ vanishes on $\{x_{1}=0\}$
and on $\bR^{d}_{+}\setminus B_{R}^{+}$.
   Then there exists  $\theta=\theta(d,\delta,K_{F},  p,   K_{0})
\in(0,1)$
such that if Assumption \ref{assump1b} ($\theta$)
is satisfied, then
$$
\int_{\bR^{d}_+}|D^{2}u|^{p}w\,dx
\leq N\int_{\bR_+^{d}}|F[u]|^{p}w\,dx
$$
\begin{equation}
                                           \label{12.15.2c}
+N \int_{\bR^{d}_+}|u|^{p}w\,dx
+N\tau_{0}^{p}\int_{\bR_+^{d}}I _{B^+_{R+R_{0} }}w\,dx,
\end{equation}
where  $N$  is a constant depending only on $d$, $\delta$, $K_F$, $K_0$,
$p$, and $R_0$.
\end{theorem}

\begin{theorem}
                             \label{theorem 4.10.4}

Take
 $K_{0}\in(1,\infty)$,   $p>   d$,
 and let $w$ be an $A_{p/ d }$-weight on $\bR^{d}_{+}$
with $[w]_{p/d}
\leq K_{0}$.
Let
\begin{equation}
                                          \label{4.11.5}
u\in   W^2_{p,w}(\bR^{d}_{+})
\end{equation}
and $u=0$ on $\{x_{1}=0\}$.
  Then there exists
  $\theta=\theta(d,\delta,K_{F},  p,   K_{0})
\in(0,1)$
such that if Assumption \ref{assump1b} ($\theta$)
is satisfied with $\tau_{0}=0$, then
\begin{equation}
                                               \label{4.11.6}
\int_{\bR^{d}_{+}}(|D^{2}u|^{p}+|Du|^{p})w\,dx
\leq N\int_{\bR^{d}_{+}}|F[u]|^{p}w\,dx
+N \int_{\bR^{d}_{+}}|u|^{p}w\,dx,
\end{equation}
and if in addition $u$ is bounded
and $F$ satisfies Assumption \ref{assumption 4.26.1} then
\begin{equation}
                                               \label{4.11.7}
\int_{\bR^{d}_{+}}(|D^{2}u|^{p}+|Du|^{p}+|u|^{p})w\,dx
\leq N\int_{\bR^{d}_{+}}|F[u]-u|^{p}w\,dx,
\end{equation}
where  the constants $N$
depend  only on $d$, $\delta$, $K_F$, $K_0$,
$p$, and $R_0$.

\end{theorem}

Proof.
Lemma
\ref{lemma 12.15.4} has a natural half space analog
and as in the case of \eqref{4.8.1} it suffices to estimate
$|D^{2}u|$.
 We prove this estimate in the same way
as in the case of
\eqref{4.8.1} by taking the same function $u_{n}$
but substituting it into \eqref{12.15.2c}
instead of \eqref{12.15.2}.

To prove \eqref{4.11.7}, it suffices to apply \eqref{4.24.6}
to
the odd extension of  $u$ and the even
 extension of $w$
across $\{x_{1}=0\}$ and use the fact that
so extended $w$ is in $A_{p}(\bR^{d } )$
with its norm controlled by
$K_{0}$.
The theorem
is proved.

The following theorem is proved in the same way as
Theorem \ref{theorem 4.10.2}, by taking into account that
$\hat x_{1}^{q}$ are $A_{p}$-weights
on $\bR^{d}_{+}$ for $q\in(-1,p-1)$.

\begin{theorem}
                                  \label{theorem 4.12.2}
Let $p_{1},p_{2}>d$, $q\in (-1,p_{2} /d -1)$, and let $u\in
C^{1,1}(\bR^{d}_{+})$
have bounded support and $u=0$ on $\{x_{1}=0\}$.
Let $F$ satisfy Assumption \ref{assumption 4.26.1}.
Then there exists
  $\theta=\theta(d,\delta, q, p_{1},p_{2}  )
\in(0,1)$
such that if Assumption \ref{assump1} ($\theta$)
is satisfied with $\tau_{0}=0$, then
 there is a constant
$N$, depending only on $d$, $\delta$, $q$, $p_{1}$, $p_{2}$,
  and $R_0$, such that
$$
\int_{0}^{\infty}\hat x_{1}^{q }
\Big(\int_{\bR^{d-1}}\Big[|  D^{2}u|+|D u|+|u|
\big]^{p_{1}}\,dx'
\Big)^{p_{2}/p_{1}}\,dx_{1}
$$
\begin{equation}
                                              \label{4.12.4}
\leq N\int_{0}^{\infty}\hat x_{1}^{q }
\Big(\int_{\bR^{d-1}}|  F[u]-u |^{p_{1}}\,dx'
\Big)^{p_{2}/p_{1}}\,dx_{1}.
\end{equation}
\end{theorem}

\begin{remark}
Estimate \eqref{4.12.4} also holds with
$x_{1}$ in place of $\hat x_{1}$. In such a situation
assume that $F(\sfu'',x)$ is independent of $x$
and is positive homogeneous of
degree one with respect to $\sfu''$. Then
scalings: $x\to c x$, immediately leads to
$$
\int_{0}^{\infty} x_{1}^{q }
\Big(\int_{\bR^{d-1}} |  D^{2}u| ^{p_{1}}\,dx'
\Big)^{p_{2}/p_{1}}\,dx_{1}
$$
$$
\leq N\int_{0}^{\infty} x_{1}^{q }
\Big(\int_{\bR^{d-1}}|  F[u]  |^{p_{1}}\,dx'
\Big)^{p_{2}/p_{1}}\,dx_{1}
$$
for any $q\in (-1,p_{2} /d -1)$ and $u\in
C^{1,1}(\bR^{d}_{+})$
with bounded support vanishing on
  $\{x_{1}=0\}$.
\end{remark}

As before, by using a localization argument, we
obtain the following estimate.

\begin{theorem}
                                      \label{theorem 12.16.1bb}
Take $x_0\in \bR^d_+$, $R\in(0,\infty)$, $K_{0}\in(1,\infty)$, $p>d$ and let
$w$ be an $A_{p/ d}$-weight with $[w]_{p/d}\leq K_{0}$. Suppose that
$D^{2}u\in L_{p,w }(B_R^+(x_0))$
and $u$ vanishes  on $\{x_{1}=0\}
\cap  B_R^+(x_0)$ if this set is nonempty.
  Then there exists
  $\theta=\theta(d,\delta,K_{F}, p,  K_{0})
\in(0,1) $
such that, if Assumption \ref{assump1bhomo} ($\theta$)
is satisfied, then for any $r\in(0,R)$
$$
\int_{B_{r}^+(x_0)}|D^{2}u|^{p}w\,dx
\leq N\int_{B_{R}^+(x_0)}|F[u]|^{p}w\,dx
$$
\begin{equation}
                                               \label{4.24.5}
+N \int_{ B_R^+(x_0)
}((R-r)^{-1}|Du|+((R-r)^{-2} +1 )|u|)^{p}w\,dx,
\end{equation}
where  $N$  is a constant depending
only on $d$, $\delta$, $K_F$, $K_0$, $p$, and $R_0$.
\end{theorem}

The proofs of the  next two theorems
 are  obtained by closely following the proof of Theorem
\ref{theorem 12.29.1} (with the lemmas proceeding it)
 with one distinction that,
since we do not have global solvability in Sobolev
spaces for equations in $B_{R}^+(x_0)$ if $B_{R}^+(x_0)
\not\subset\bR^{d}_{+}$, we take a smooth subdomain
of $B_{R}^+(x_0)$
containing $B_{r}^+(x_0)$ and conduct the corresponding
argument in the proof of   Theorem
\ref{theorem 12.29.1} with this subdomain in place of
$B_{R}$.

\begin{theorem}
                                    \label{theorem 3.25.1}
 Take $R\in(0,\infty)$, $r\in (0,R)$, $p>  d$, $p_{1},p_{2}>d$, and
 $u\in W^{2}_{p}(B_{R}^+ )$.
Suppose that $u$ vanishes  on $\{x_{1}=0\}$.
Finally,
take $q\in(-1,p_{2} /d -1)$ and let $F$ satisfy Assumption \ref{assumption 4.26.1}.
 Then there exists
  $\theta=\theta(d,\delta, p,  q,p_{1})
\in(0,1) $
such that, if Assumption  \ref{assump1bhomo}  ($\theta$)
is satisfied, then
\begin{multline*}
\int_{0}^{\infty}x_{1}^{q}\Big(\int_{\bR^{d-1}}
I_{B_{r}^{+}}(|D^{2}u|^{p_{1}}+|Du|^{p_{1}})\,dx'\Big)^{p_{2}/p_{1}}dx_{1}
\\
\leq N\int_{0}^{\infty}x_{1}^{q}\Big(\int_{\bR^{d-1}}
I_{B_{R}^{+}}|F[u]|^{p_{1}}\,dx'\Big)^{p_{2}/p_{1}}dx_{1}
\\
+N\int_{0}^{\infty}x_{1}^{q}\Big(\int_{\bR^{d-1}}
I_{B_{R}^{+}}|u|^{p_{1}}\,dx'\Big)^{p_{2}/p_{1}}dx_{1},
\end{multline*}
where the constants
 $N$ depend only on  $r$, $R$,
$d$, $\delta$, $p$, $p_{1}$, $p_{2}$, $q$, and $R_0$.
\end{theorem}

\begin{theorem}
                                   \label{theorem 12.29.1bdry}

Take $x_0\in   \bR^d_+$, $R\in(0,\infty)$, $r\in (0,R)$, $p>  d$ and
take  $p_{i}>d$,  $i= 1, 2,\ldots,d$.
Assume that $u\in W^{2}_{p}(B_{R}^+(x_0))$ and $u$ vanishes  on $\{x_{1}=0\}$.
Then there exists
  $\theta=\theta(d,\delta, p,p_{1},\ldots,p_{d})
\in(0,1) $
such that, if Assumption  \ref{assump1bhomo}  ($\theta$)
is satisfied and Assumption \ref{assumption 4.26.1}
is satisfied as well, then
$$
\| I_{B_{r}^+(x_0)} D^{2}u,I_{B_{r}^+(x_0)} Du  \|_{L_{p_{1},\ldots,p_{d}}(\bR^{d})}
\leq N\|I_{B_{R}^+(x_0)} F[u] \,\|_{L_{p_{1},\ldots,p_{d}}(\bR^{d})}
$$
$$
+N\|I_{B_{R}^+(x_0)} u
\|_{L_{p_{1},\ldots,p_{d}}(\bR^{d})},
$$
where the constants
 $N$ depend only on  $r$, $R$,
$d$, $\delta$, $p_1,\ldots,p_d$, and $R_0$.
\end{theorem}

\mysection{Parabolic case}
                                   \label{section 4.23.1}

We concentrate our attention here  on
$$
\bR^{d+1}_{+}=\{(t,x):t\geq0,x\in\bR^{d}\},
$$
and on functions defined on it.

For $(t,x)\in \bR^{d+1}_{+} $ introduce
$$
C_{r}(t,x)=[t,t+r^{2})\times B_{r}(x),\quad C_{r}=C_{r}(0,0).
$$
We consider a function $F(\sfu'',t,x)$, $\sfu''\in\bS$,
$(t,x)\in\bR^{d+1}_{+} $,
on which we will impose  some of the following assumptions.

\begin{assumption}[$\theta$]
                             \label{assumption 12.26.1}
Assumption \ref{assump1} ($\theta$) is satisfied if we replace there $x$, $\bR^{d},B_{r}(z)$
 with $(t,x)$, $\bR^{d+1}_{+} $, $C_{r}(z)$, respectively.
\end{assumption}

 \begin{assumption}
                                         \label{assumption 4.26.1p}
Assumption \ref{assumption 4.26.1} is satisfied if we replace there $F(\cdot,x)$ with $F(\cdot,t,x)$.
\end{assumption}

 \begin{assumption}[$\theta$]
                                       \label{assumption 12.26.1homo}
Assumption \ref{assump1homo} ($\theta$) is satisfied if we replace there $x$, $\bR^{d},B_{r}(z)$
by $(t,x)$, $\bR^{d+1}_{+} $, $C_{r}(z)$, respectively.
\end{assumption}

By using similar natural substitutions, we introduce
$$
h^{\sharp}_{ \gamma,\rho}(t,x),\quad  \bM_\rho h (t,x),\quad  \text{and}\quad\bM h(t,x)
$$
on $\bR^{d+1}_{+}$
(taking only $C_{r}( t_0, x_{0})
\subset \bR^{d+1}_{+} $, $ (t_0 ,x_{0})\in \bR^{d+1}_{+} $).

Here we set $\Omega=\bR^{d+1}_{+} $ and
for $n\in\bZ$ we take $\bC_{n}$
as the collection of $(t,x)+[0,4^{-n})\times[0,2^{-n})^{d}$,
$t\in 4^{-n}\{0,1,\ldots\}$, $x\in 2^{-n}\bZ^{d}$.
We also set $\mu$ to be Lebesgue measure
on $\bR^{d+1}_{+}$ and $\bL$ to be the set of continuous
functions on $\bR^{d+1}_{+}$ with compact support.
Then
observe that relations \eqref{12.15.5}
hold again for a constant $c=c(d)\in(1,\infty)$.

In what follows in this section by $A_{p}$-weights we mean
weights on $\bR^{d+1}_{+}$ relative to the parabolic distance.

The following analog of Lemma \ref{lemma 12.15.1}
is an obvious corollary of
Lemma 3.3 of \cite{Kr_18}.

\begin{lemma}
                                        \label{lemma 12.27.1}
Let $u\in W^{1,2}_{d+1,\loc }(
\bR^{d+1}_{+})$,
$\mu\in(0,\infty)$, $\nu\geq2$,
$\xi\in(1,\infty)$.  Then there exists
$\theta=\theta(d,\delta,K_{F},\mu,\xi )\in(0,1)$ such that,
if Assumption \ref{assumption 12.26.1} $(\theta)$ is
satisfied, then one can find
$ \gamma_0=\gamma_0(d,\delta) \in(0,1)$,
$ \alpha=\alpha(d,\delta)\in(0,1)$,
such that, for $\gamma\in(0,\gamma_{0}]$, $h=D^{2}u$, and $\rho=R_{0}/\nu$,
   we have
$$ h^{\sharp}_{ \gamma,\rho}\leq
N\nu^{(d+2)/\gamma}\bM^{1/(d+1)}\big[|\partial_{t}u+F[u]|^{d+1}\big]
+N\tau_{0}\nu^{(d+2)/\gamma}
$$
\begin{equation}
                                            \label{12.27.2}
+N(\mu\nu^{(d+2)/\gamma}+\nu^{-\alpha})
\bM^{1/(\xi'(d+1))}\big[|h|^{\xi'(d+1)}\big],
\end{equation}
where $\xi'=\xi/(\xi-1) $ and the constants $N$
depend only on $d$, $K_{F}$,   and $\delta$.

\end{lemma}

Here is a parabolic analog of Lemma
\ref{lemma 12.15.4}.

\begin{lemma}
                           \label{lemma 12.27.2}
(i) There exists a constant
$\gamma_{0}=\gamma_{0}(d,\delta)\in(0,1)$
such that for any $\gamma\in(0,\gamma_{0}]$,
$\rho>0$, and $u\in W^{1,2}_{d+1,\loc}(\bR^{d+1}_{+})$, we have
$$
\bM^{1/\gamma}_{\rho}(|D^{2}u|^{\gamma})\leq N \bM^{1/(d+1)}_{\rho}(|
\partial_{t}u+F[u]|^{d+1})
$$
\begin{equation}
                                     \label{12.27.3}
+N   \rho^{-1} \bM ^{1/(d+1)}_{ \rho}(|Du|^{d+1})
+N \rho^{-2} \bM^{1/(d+1)}_{\rho}(|u|^{d+1}),
\end{equation}
where the constants $N$ depend only on
$d$, $\delta$, and $K_F$.

(ii)  For any $\rho>0$, $p\in [1,\infty)$,
and $u\in W^{1,2}_{p,\text{loc}}
(\bR^{d+1}_{+})$, we have
\begin{equation}
                                            \label{12.27.4}
  \bM _{\rho}(|Du|^{ p })\leq N\bM^{1/2}_{\rho}(|D^{2}u|^{p})
\bM^{1/2}_{\rho}(| u|^{p})+  N\rho^{-p} \bM _{\rho}(| u|^{p}),
\end{equation}
where the constants $N$ depend only on $d$ and $p$.

(iii) For any $\rho>0$, $K_{0}$,
$p\in (1,\infty)$, $w\in A_p$ with $[w]_{p}\le K_0$, and $u\in W^{1,2}_{p,w}
(\bR^{d+1}_{+})$, we have
\begin{equation}
                                            \label{1.38p}
\int_{\bR^{d+1}_{+}}|Du|^{p}w\,dxdt\leq \rho^{p}
\int_{\bR^{d+1}_{+}}
\big |D^{2}u|^{p}w\,dxdt+N\rho^{-p}
\int_{\bR^{d+1}_{+}}
\big |u|^{p}w\,dxdt,
\end{equation}
where $N$ depends only on $d$, $p$, and $K_0$.
\end{lemma}

Proof. First write $F[u]=a^{ij}D_{ij}u$
and take $r\geq\rho$
  and a function $\zeta\in C^{\infty}_{0}(\bR^{d+1})$
such that $\zeta=1$ on $C_{r}$, $\zeta=0$ on
$\partial'C_{2r}$, and
 $$
|D\zeta|\le N/r\le N/\rho,\quad |\zeta_t|+|D^2\zeta|\le N/r^2\le N/\rho^2.
$$

Then
by Lemma 5.5 of \cite{DKL_12}
$$
\dashint_{C_{r}}|D^{2}u|^{\gamma}\,dxdt
\leq \dashint_{C_{2r}}|D^{2}(\zeta  u)|^{\gamma}\,dxdt
$$
$$
\leq N\Big( \dashint_{C_{2r}}|\partial_{t}(\zeta u)
+\zeta  F[u]+a^{ij}2D_{i}\zeta
D_{j}u+ ua^{ij}D_{ij}\zeta )|^{d+1}\,dx dt \Big)^{\gamma/(d+1)}.
$$
 The rest is identical to the proof of
Lemma \ref{lemma 12.15.4}. The lemma is proved.

\begin{theorem}
                                      \label{theorem 12.29.1b}
Take $R\in(0,\infty)$, $K_{0}\in(1,\infty)$,
$p>d+1$
and
let $w$ be an $A_{p/ (d+1) }$-weight
with
$[w]_{p/ (d+1)}
\leq K_{0}$.
Suppose that $D^{2}u \in L_{p,w}(\bR^{d+1}_{+} )$, and that $u$
vanishes in
$\bR^{d+1}_{+} \setminus C_{R}$.   Then there exists
  $\theta=\theta(d,\delta,K_{F},   p,  K_{0})
\in(0,1) $
such that, if Assumption \ref{assumption 12.26.1} ($\theta$)
is satisfied, then
$$
\int_{\bR^{d+1}_{+} }|D^{2}u|^{p}w\,dxdt
\leq N\int_{\bR^{d+1}_{+} }|\partial_{t}u+F[u]|^{p}w\,dxdt
$$
\begin{equation}
                                                    \label{12.29.2}
+N \int_{\bR^{d+1}_{+} }|u|^{p}w\,dxdt
+N\tau_{0}^{p}\int_{\bR^{d+1}_{+} }I _{C_{R+R_{0} }}w\,dxdt,
\end{equation}
where  $N$  depends only on $d$, $\delta$,
 $K_F$, $K_0$, $p$, and $R_0$.
\end{theorem}

The proof of this theorem is practically the same as
that of Theorem \ref{theorem 12.15.1}.

To prove a parabolic analog of Theorem \ref{theorem 4.8.1}
we need the following analog of Lemma \ref{lemma 4.8.1}.

\begin{lemma}
                                         \label{lemma 3.30.1}
Let $u\in W^{1,2}_{d+1,\loc}(\bR^{d+1}_{+})$
be a bounded function
and $a=(a^{ij}(t,x))$ be an $\bS_{\delta}$-valued
function on $\bR^{d+1}$. Also let $p>d+1$.
Then
$$
|u(0)|\leq N\big(\bM(|\partial_{t}u+a^{ij}D_{ij}u
-u|^{p})(0)\big)^{1/p},
$$
where $N=N(d,\delta,p)$.
\end{lemma}

Proof. Let $G(s,t,x,y)$ be a  Green's function of
$L:=\partial_{t} +a^{ij}D_{ij}-1$ in $\bR^{d+1}_+$
and introduce $f=-Lu$. Then
for $G(t,y):= G(0,t,0,y) $ we have
$$
u(0)=\int_{0}^{\infty}\int_{\bR^{d}}
G( t,y)f(t,y)\,dydt.
$$
Hence,
$$
|u(0)|\leq \int_{0}^{\infty}\int_{\bR^{d}}
G( t,y)|f(t,y)|\,dydt.
$$

We are going to use the following estimate
easily obtained, say, by probabilistic arguments:
for any $\alpha\geq0$
\begin{equation}
                                          \label{3.31.1}
\int_{\bR^{d+1}_{+}}G(t,y)(t^{2\alpha}+|y|^{\alpha})
\,dydt\leq N(\alpha,d,\delta).
\end{equation}

Observe that for any $h(t,y)\geq0$ and $\alpha>0$
$$
\int_{1}^{\infty}r^{-\alpha-1}\Big(\int_{C_{r}}h (t^{2\alpha}\vee
|y|^{\alpha}\vee1)\,dydt\Big)dr
$$
$$
=\int_{\bR^{d+1}_{+}}h(t^{2\alpha}\vee
|y|^{\alpha}\vee 1)\Big(\int_{t^{2}\vee|y|\vee 1 }^{\infty}
r^{-\alpha-1}\,dr\Big)
dydt
=\frac{1}{\alpha}\int_{\bR^{d+1}_{+}}h \,
dydt.
$$
Furthermore, by using \eqref{3.31.1}, H\"older's inequality,
and the parabolic Aleksandrov estimate, for  $q=p/(d+1)>1$,
we get
$$
\int_{C_{r}}G(t,y)|f(t,y)|(t^{2\alpha}\vee
|y|^{\alpha}\vee 1)\,dydt\leq N\Big(\int_{C_{r}}G(t,y)|f(t,y)|^{q}
\,dydt\Big)^{1/q}
$$
$$
\leq N\Big(\int_{C_{r}}|f(t,y)|^{p}
\,dydt\Big)^{1/p}\leq Nr^{(d+2)/p}\big(\bM(|f|^{p})(0)\big)^{1/p}.
$$
For $(d+2)/p-\alpha-1<-1$, we get the desired result by
integrating in $r\in [1,\infty)$ and collecting the above
estimates. The lemma is proved.

Now by combining Theorem \ref{theorem 12.29.1b}
and Lemmas \ref{lemma 3.30.1} and \ref{lemma 12.27.2}
we get  the following in the same way as
Theorem \ref{theorem 4.8.1}.

\begin{theorem}
                                    \label{theorem 4.24.1}
Let $\tau_{0}=0$ and take
 $K_{0}\in(1,\infty)$.
Let   $p>   d+1$
 and let $w$ be an $A_{p/ (d+1) }$-weight
with $[w]_{p/(d+1)}
\leq K_{0}$.
Let
\begin{equation}
                                          \label{4.24.2}
u\in   W^{1,2}_{p,w }(\bR^{d+1}_{+}).
\end{equation}

  Then there exists
  $\theta=\theta(d,\delta,K_{F},  p,   K_{0})
\in(0,1)$
such that if Assumption \ref{assumption 12.26.1} ($\theta$)
is satisfied, then
$$
\int_{\bR^{d+1}_{+}}(|D^{2}u|^{p}+|Du|^{p})w\,dxdt
$$
\begin{equation}
                                               \label{4.24.3}
\leq N\int_{\bR^{d+1}_{+}}|\partial_{t}u+F[u]|^{p}w\,dxdt
+N \int_{\bR^{d+1}_{+}}|u|^{p}w\,dxdt,
\end{equation}
and if, in addition, $u$ is bounded and Assumption
\ref{assumption 4.26.1p}
is satisfied then
\begin{equation}
                                               \label{4.24.4}
\int_{\bR^{d+1}_{+}}(|D^{2}u|^{p}+|Du|^{p}+|u|^{p})w\,dxdt
\leq N\int_{\bR^{d+1}_{+}}|\partial_{t}u+F[u]-u|^{p}w\,dxdt,
\end{equation}
where  the constants $N$
depend  only on $d$, $\delta$, $K_F$, $K_0$,
$p$, and $R_0$.

\end{theorem}

To state an analog of Theorem \ref{theorem 4.9.2}
order the set of coordinates
$(t,x)=(t,x_{1},\ldots,x_{d})$ arbitrarily
as $(\tilde x_{0}, \ldots,\tilde x_{d})$. Then we have the following result.

\begin{theorem}
                                    \label{theorem 4.24.3}
Let  $\tau_{0}=0$  and
take  $p_{i}>d+1$,  $i= 0, 1,\ldots,d$.
Assume that $u\in W^{1,2}_{1,\loc}(\bR^{d+1}_{+})$
  and
Assumption \ref{assumption 4.26.1p} is satisfied.
Then there exists $\theta=\theta(d,\delta,    d
,p_{0},\ldots,p_{d} )\in(0,1)$ such that, if Assumption
\ref{assumption 12.26.1} ($\theta$)
is satisfied, then
\begin{equation}
                                      \label{4.9.5p}
\|   D^{2}u,Du,u \|_{L_{p_{0},\ldots,p_{d}}(\bR^{d+1}_{+})}
\leq N\|\partial_{t}u+ F[u]-u \,\|_{L_{p_{0},\ldots,p_{d}}
(\bR^{d+1}_{+})}
\end{equation}
provided  that the left-hand side is finite,
where
$$
\| f
\|_{L_{p_{0},\ldots,p_{d}}(\bR^{d+1}_{+})}^{p_{d}}
$$
\begin{equation}
                                                    \label{3.28.1}
:=\int_{\bR  }\Big(\cdots
\Big(\int_{\bR }
\Big(\int_{\bR } |fI_{\bR^{d+1}_{+}}|^{p_{0}}\,d\tilde  x_{0}
\Big)^{p_{1} /p_0 }\,d\tilde x_{1}\Big)^{p_{2}/p_{1}}
\cdots\Big)^{p_{d} /p_{d-1} }\,d\tilde  x_{d},
\end{equation}
and the constant  $N$ depends only on $d$, $\delta$,   $p_1,\ldots,p_d$, and $R_0$.
\end{theorem}

One proves this result in the same way as Theorem \ref{theorem 4.9.2}
taking care of defining the mollified functions
$u(t,x)$ by averaging the values of $u$ with higher values
of $t$ in order not to bother about the fact that $u$
may not be defined for negative $t$.

Then one derives an obvious analogs
 of Lemmas \ref{lemma 12.16.1}
and \ref{lemma 4.10.6}
and, by using Theorem 1.9 of \cite{Kr_18}
(see also Remark 1.11 there)
in place of  Theorem 2.1 of \cite{Kr_13},
 one arrives at an analog of
 Theorem \ref{theorem 4.11.1}.

\begin{theorem}
                                    \label{theorem 4.24.4}
Theorem \ref{theorem 4.24.1} remains true
if condition \eqref{4.24.2} is replaced with
$u\in \bigcap_{R>0}  W^{1,2}_{p,w}(C_{R})$ provided,
 additionally,
that $F$ is positive homogeneous
of degree one with respect to $\sfu''$.
In particular, if $u$ is bounded and
 the right-hand side of
\eqref{4.24.4} is finite, then
$u\in W^{1,2}_{p,w}(\bR^{d+1}_{+})$.

\end{theorem}
\begin{remark}
                                      \label{remark 4.24.1}

Generally,  \eqref{4.24.4} may fail if $u$ is unbounded.
Indeed, if $d=1$ and $F[u]=u''$, the function $e^{x}$
satisfies $\partial_{t}u+F[u]-u=0$ and is nonzero.
\end{remark}

Then from an analog of Lemma \ref{lemma 4.10.6}
one derives the following analog of
Theorem \ref{theorem 12.29.1}. The only difference
in the proofs worth noting is that
one should use the existence Theorem 1.9 of \cite{Kr_18}
in place of Theorem
2.1 of \cite{Kr_13}.

\begin{theorem}
                                   \label{theorem 12.29.1p}

Take $R\in(0,\infty)$,  $r\in (0,R)$,
$p>d+1$, and $p_i>d+1$ for $i=0,1,\ldots,d$. Assume that
$u\in W^{1,2}_{p}(C_{R})$.
 Then there exists
  $\theta=\theta(d,\delta,  p ,p_{0},\ldots,p_{d})
\in(0,1) $
such that, if Assumptions  \ref{assumption 12.26.1homo}
 ($\theta$)  and \ref{assumption 4.26.1p}
are satisfied, then
$$
\| I_{C_{r}}  D^{2}u,I_{C_{r}}
 D u \|_{L_{p_{0},\ldots,p_{d}}(\bR^{d+1}_{+})}
\leq N\|I_{C_{R} } F[u] \,\|_{L_{p_{0},\ldots,p_{d}}
(\bR^{d+1}_{+})}
$$
\begin{equation}
                                      \label{2.20.3p}
+N\|I_{C_{R} } u
\|_{L_{p_{0},\ldots,p_{d}}(\bR^{d+1}_+)},
\end{equation}
where  the constants $N$ depend only on
$r$, $R$, $d$, $\delta$,  $p_0,\ldots,p_d$,
and $R_0$.
\end{theorem}

\mysection{Parabolic case in a half-space}
                                \label{sec7}

Here we consider functions on
$$
\bR^{d+1}_{+,+}=\{(t,x):t\geq 0,x_{1}\geq0, x'\in\bR^{d-1}\}.
$$

We concentrate on parabolic equations in
  $\bR^{d+1}_{+,+}$
 with
 zero Dirichlet
boundary condition  and prove boundary estimates with $A_{p}$-weights.

For $(t,x)\in \bR^{d+1}_{+,+}$ and $r>0$ denote
$$
C^+_{r}(t,x)=[t,t+r^{2})\times B^+_{r}(x),\quad C^+_{r}=C^+_{r}(0,0).
$$
and consider a function $F(\sfu'',t,x)$ given for $(t,x)\in
\bR^{d+1}_{+,+}$
and $ \sfu'' \in\bS$.

We use the following assumptions.

\begin{assumption}[$\theta$]
                                       \label{assumption 3.27.1}
Assumption \ref{assump1} ($\theta$) is satisfied if we replace
there $x$, $\bR^{d},B_{r}(z)$  with $(t,x)$, $\bR^{d+1}_{+,+} $,
$C^{+}_{r}(z)$, respectively.
\end{assumption}

 \begin{assumption}[$\theta$]
                                       \label{assumption 3.27.1homo}
Assumption \ref{assump1homo} ($\theta$) is satisfied if we replace
there $x$, $\bR^{d},B_{r}(z)$ with $(t,x)$, $\bR^{d+1}_{+,+} $,
$C^{+}_{r}(z)$, respectively.
\end{assumption}

Accordingly, we introduce
$$
h^{\sharp}_{ \gamma,\rho}(t,x),
\quad \bM_\rho h(t,x),\quad  \text{and}\quad\bM h(t,x)
$$
on $\bR^{d+1}_{+,+}$ (by taking $C_{r}^+(t_{0},x_{0})
\subset \bR^{d+1}_{+,+}$, $(t_{0},x_{0})\in  \bR^{d+1}_{+,+}
$). Here the underlying set $\Omega$ is taken to be
$\bR^{d+1}_{+,+}$ and  the $\bC_{n}$'s are
 the parts of the $\bC_{n}$'s from the
beginning of Section \ref{section 4.23.1} which belong
to that $\Omega$.

In what follows by $A_{p}$-weights we mean
weights on $\bR^{d+1}_{+,+}$ relative to the parabolic distance.

From Lemma 4.1 of \cite{Kr_18} and the proof of Lemma 3.3 of
\cite{Kr_18},
 we can easily obtain a boundary analog of Lemma \ref{lemma
12.15.1}. This together with a boundary analog of Lemma
\ref{lemma 12.15.4},
by relying on Corollary \ref{corollary 12.3.1}, gives the following boundary estimate
corresponding to Theorems \ref{theorem 12.15.1}
and \ref{theorem 12.29.1b}.

\begin{theorem}
                                      \label{theorem 12.15.1c}

Take $R\in(0,\infty)$,  $K_{0}\in(1,\infty)$. Let $p>d+1$  and let
$w$ be an $A_{p/(d+1)}$-weight on $\bR^{d+1}_{+,+}$ with $[w]_{p/(d+1)}
\leq K_{0}$. Suppose that
$$
D^{2}u \in L_{p,w}(\bR_{+,+}^{d+1})
$$   and $u$ vanishes on $\{x_{1}=0\}$ and on $
\bR^{d+1}_{+,+}\setminus C_{R}^{+}$.
   Then there exists  $\theta=\theta(d,\delta,K_{F},  p,   K_{0})
\in(0,1)$
such that if Assumption \ref{assumption 3.27.1} ($\theta$)
is satisfied, then
$$
\int_{\bR^{d+1}_{+,+}}|D^{2}u|^{p}w\,dxdt
\leq N\int_{\bR_{+,+}^{d+1}}|
\partial_{t}u+F[u]|^{p}w\,dxdt
$$
\begin{equation}
                                                    \label{12.15.2cc}
+N \int_{\bR^{d+1}_{+,+}}|u|^{p}w\,dxdt
+N\tau_{0}^{p}\int_{\bR_{+,+}^{d+1}}I _{C^+_{R+R_{0} }}w\,dxdt,
\end{equation}
where  $N$  is a constant depending only on $d$, $\delta$, $K_F$, $K_0$,
$p$, and $R_0$.
\end{theorem}

 By taking into account what was said before
Theorems \ref{theorem 3.25.1} and \ref{theorem 12.29.1bdry}
and using the solvability of $\partial_{t}u+F[u]=f$ in smooth cylinders
(see Theorem 1.9 and Remark 1.11 of \cite{Kr_18}), we
have the following boundary estimates in mixed-norm spaces.

\begin{theorem}
                                    \label{theorem 3.28.2}

Take $R\in(0,\infty)$, $r\in (0,R)$, $p>  d+1$, $p_{1},p_{2}>d+1$, and
 $u\in W^{ 1, 2}_{p}(C_{R}^+ )$.
Suppose that $u$ vanishes  on $\{x_{1}=0\}$.
Finally,
take $q\in(-1,p_{1} /(d+1) -1)$ and let $F$ satisfy Assumption \ref{assumption 4.26.1p}.
 Then there exists
  $\theta=\theta(d,\delta, p,  q)
\in(0,1) $
such that, if Assumption  \ref{assumption 3.27.1homo}  ($\theta$)
is satisfied, then
\begin{multline}
                                               \label{3.29.1}
\int_{0}^{\infty}\Big(\int_{\bR^{d}_{+}}
I_{C_{r}^{+}}x_{1}^{q}(|D^{2}u|^{p_{1}} +|Du|^{p_{1}})\,dx\Big)^{p_{2}/p_{1}}dt
\\
\leq N\int_{0}^{\infty}\Big(\int_{\bR^{d}_{+}}
I_{C_{R}^{+}}x_{1}^{q}|\partial_{t}u+F[u]|^{p_{1}}\,dx\Big)^{p_{2}/p_{1}}
dt
\\
+N\int_{0}^{\infty} \Big(\int_{  \bR^{d}_{+}}
I_{C_{R}^{+}} x_{1}^{q}|u|^{p_{1}}\,dx\Big)^{p_{2}/p_{1}}dt,
\end{multline}
\begin{multline}
                                               \label{3.29.2}
\int_{\bR^{d}_{+}}x_{1}^{q}\Big(\int_{0}^{\infty}
I_{C_{r}^{+}}(|D^{2}u|^{p_{2}} +|Du|^{p_{2}})\,dt\Big)^{p_{1}/p_{2}}dx
\\
\leq N\int_{\bR^{d}_{+}}x_{1}^{q}\Big(\int_{0}^{\infty}
I_{C_{R}^{+}} |\partial_{t}u+F[u]|^{p_{2}}\,dt\Big)^{p_{1}/p_{2}}
dx
\\
+N\int_{\bR^{d}_{+}}x_{1}^{q}\Big(\int_{0}^{\infty}
I_{C_{R}^{+}}  |u|^{p_{2}}\,dt\Big)^{p_{1}/p_{2}}dt,
\end{multline}
where  the constants
 $N$ depend only on  $r$, $R$,
$d$, $\delta$, $p$, $p_{1}$, $p_{2}$, $q$,  and $R_0$.
\end{theorem}

\begin{theorem}
                                   \label{theorem 3.13bdry}

Take $(t_0,x_0)\in  \bR^{d+1}_{+,+}$, $R\in(0,\infty)$, $r\in (0,R)$,
 $p>  d+1$ and
take  $p_{i}>d+1$,  $i= 0, 1,\ldots,d$.
Assume that $u\in  W^{1,2}_{p} (C^+_{R}(t_0,x_0))$
and $u$ vanishes on $\{x_{1}=0\}$.
Then there exists
 $\theta=\theta(d,\delta, p , p_0
 ,\ldots,p_d  )
\in(0,1)$
such that if Assumption  \ref{assumption 3.27.1homo}  ($\theta$)
is satisfied, then (the mixed norms below are taken from
\eqref{3.28.1})
$$
\| I_{C_{r}^+(t_0,x_0)}  D^{2}u \|_{L_{p_{0},\ldots,p_{d}}(\bR^{d+1})}
\leq N\|I_{C_{R}^+(t_0,x_0)} F[u] \,\|_{L_{p_{0},\ldots,p_{d}}(\bR^{d+1})}
$$
$$
+N\|I_{C_{R}^+(t_0,x_0)} u
\|_{L_{p_{0},\ldots,p_{d}}(\bR^{d+1})},
$$
where the constants
 $N$ depend only on  $r$, $R$,
$d$, $\delta$,  $p$,  $p_0,\ldots,p_d$, and $R_0$.
\end{theorem}

 To further estimate the lower-order terms on the right-hand
sides of the estimates above, we need the following fact.

By using the odd extension of $u$ and the even extension of $w$
across $\{x_{1}=0\}$ and using the fact that
so extended $w$ is in $A_{p}(\bR^{d+1}_{+})$ with its norm controlled by
$K_{0}$,
from Lemmas \ref{lemma 3.30.1}
we get the following corollary in which
$$
 W^{1,2}_{d+1 ,\loc  }(\bR_{+,+}^{d+1})=
\bigcap_{R>0} W^{1,2}_{d+1 }(C_{R}^{+}).
$$

\begin{corollary}
                           \label{corollary 4.1.1}
 Let $K_0\in (1,\infty)$,
$p>d+1$ and let $u\in W^{1,2}_{d+1,\loc}(\bR^{d+1}_{+,+})$ be a bounded
function and $a$ be an $\bS_{\delta}$-valued
function on $\bR^{d+1}_{+,+}$.
Let $w\in A_{p/(d+1)}$  on $\bR_{+,+}^{d+1}$ with
$[w]_{p/(d+1)}\le K_0$  and let $u=0$ for
 $x_{1}=0$. Then
$$
\int_{\bR^{d+1}_{+,+}}|u|^{p} w \,dxdt
\leq
N\int_{\bR^{d+1}_{+,+}}|\partial_{t}u
+a^{ij}D_{ij}u-u|^{p} w \,dxdt,
$$
where $N=N(d,\delta,p,K_0)$.
\end{corollary}

We are now ready to prove the following theorem.
\begin{theorem}
                                  \label{theorem 4.1.1}

Let $K_0\in (1,\infty)$,   $p>d+1$,
 $w\in A_{p /(d+1) }$  on $\bR_{+,+}^{d+1}$ with
$[w]_{p/(d+1)}\le K_0$, and $u\in
W^{1,2}_{p,w}(\bR_{+,+}^{d+1})$ vanishing on $\{x_1=0\}$.
Let $F$ satisfy Assumption \ref{assumption 4.26.1p}.
Then  there exists  $\theta=\theta(d,\delta, p, K_0)
\in(0,1)$
such that if Assumption  \ref{assumption 3.27.1homo}  ($\theta$)
is satisfied, then
$$
\int_{\bR_{+,+}^{d+1}}\big(|D^{2}u|^{p}
+|Du|^{p}+|u|^{p}\big)w\,dxdt
\leq NI,
$$
where
$$
I=\int_{\bR_{+,+}^{d+1}}
|\partial_{t}u+F[u]-u|^{p} w\,dxdt
$$
and $N$ depends only on
$d$, $\delta$,  $K_0$, $p$, and $R_0$.
\end{theorem}

Proof. Observe that the following is a parabolic analog
of \eqref{4.24.5} for $\bR^{d}_{+,+}$:
$$
\int_{C^{+}_{1}(t_{0},x_{0})}|D^{2}u|^{p}w\,dxdt
\leq N\int_{C^{+}_{2}(t_{0},x_{0})}|\partial_{t}u+F[u]|^{p}w\,dxdt
$$
\begin{equation}
                                              \label{4.1.2}
+N \int_{C^{+}_{2}(t_{0},x_{0})} (|Du|+ |u|)^{p}w\,dxdt.
\end{equation}
The way to obtain it from
Theorem \ref{theorem 12.15.1c}
is described in the proof of Theorem \ref{theorem 12.16.1bb}
and could be easily mimicked in the parabolic setting.

By integrating both sides of \eqref{4.1.2} with respect to $(t_{0},x_{0})
\in\bR^{d+1}_{+,+}$ we get
$$
\int_{\bR^{d+1}_{+,+}}|D^{2}u|^{p}w\,dxdt
\leq N\int_{\bR^{d+1}_{+,+}}|\partial_{t}u+F[u]|^{p}w\,dxdt
$$
$$
+N \int_{\bR^{d+1}_{+,+}} (|Du|+ |u|)^{p}w\,dxdt
$$
$$
\leq NI+N \int_{\bR^{d+1}_{+,+}} (|Du|+ |u|)^{p}w\,dxdt.
$$
 By using Corollary \ref{corollary 4.1.1} and a boundary parabolic analog of
Lemma \ref{lemma 12.15.4} (iii), we arrive at
$$
\int_{\bR_{+,+}^{d+1}}\big(|D^{2}u|^{p}
+|Du|^{p}+|u|^{p}\big)w\,dxdt
\leq N \rho^{-p}
I+N\rho^{p}\int_{\bR_{+,+}^{d+1}}|D^{2}u|^{p} w \,dxdt
$$
 for any $\rho\in (0,1)$. The desired estimate follows by
taking $\rho$ sufficiently small.
The theorem is proved.

Theorems \ref{theorem 4.1.1} and \ref{theorem 12.2.1}
and the way Theorem \ref{theorem 4.9.2} is derived
immediately lead to the following.

\begin{theorem}
                                   \label{theorem 5.1.1}
 Let $ p_{1},p_{2},p_{3}>d+1$, and
 $u\in W^{ 1, 2}_{1,\loc}(\bR^{d+1}_{+,+})$.
Suppose that $u$ vanishes  on $\{x_{1}=0\}$.
Finally,
take $q\in(-1,p_{1} /(d+1) -1)$
and let $F$ satisfy Assumption \ref{assumption 4.26.1p}.
 Then there exists
  $\theta=\theta(d,\delta, p_{1},p_{2},p_{3},  q)
\in(0,1) $
such that, if Assumption  \ref{assumption 3.27.1homo}  ($\theta$)
is satisfied, then
$$
\int_{0}^{\infty}\Big(\int_{\bR^{d-1}}\Big(\int_{0}^{\infty}
 x_{1}^{q}\big[|D^{2}u|+|D u| +| u|
\big]^{p_{1}}\,dx_{1}\Big)^{p_{2}/p_{1}}
dx'\Big)^{p_{3}/p_{2}}dt
$$
\begin{equation}
                                               \label{5.1.1}
\le N\int_{0}^{\infty}\Big(\int_{\bR^{d-1}}\Big(\int_{0}^{\infty}
 x_{1}^{q}|\partial_{t}u+F[u]-u|^{p_{1}}\,dx_{1}\Big)^{p_{2}/p_{1}}
dx'\Big)^{p_{3}/p_{2}}dt,
\end{equation}
provided that the left-hand side is finite,
where $N$ depends only on
$d$, $\delta$, $ p_{1}$, $p_{2}$, $p_{3}$, $  q$,  and $R_{0}$.

\end{theorem}

The one-dimensional
 example of $F[u]=D^{2}u$ and $u(t,x)=\sinh x$ shows that \eqref{5.1.1}
is wrong without the additional assumption on its
left-hand side.

\begin{remark}
The reader understands that one has similar estimates
for the integrals with respect to $x_{1}$, $x'$, and $t$
mixed in any other order.
\end{remark}

\begin{remark}
In \cite{NK_09} the authors consider linear $F$ with  coefficients
depending only on time in a measurable way and prove
a priori estimates similar to the one in Theorem \ref{theorem 5.1.1},
however, for any $p_{1}=p_{2}, p_{3}>1$ and $q\in(-1,2p_{1}-1)$.
 The latter range
is much wider than ours $(-1,p_{1} /(d+1) -1)$,
 but our operators are much more general
and we have three   integrals.

It is worth noting that the range $(p_{1}-1,2p_{1}-1)$ was used
in \cite{Kr_01} to   build the solvability theory
of parabolic equations  in Sobolev spaces with weights
with the highest order of derivatives being
an arbitrary given number:
 positive, negative,
integral or fractional.

\end{remark}

\mysection{Appendix}

Here we take
$\Omega=\Omega^{1}\times\cdots\times\Omega^{d}$,
 where $\Omega^j=\bR$ or $\bR_+$,
$j=1,\ldots,d$   and
let $\mu$ to be the Lebesgue measure on $\Omega$.
We take integers $0= l_{0}<l_{1}<\ldots<l_{m}=d$ and express points in $\Omega$ as
$$
x=(x_{1},\ldots,x_{d})=(\check x_{1},\ldots,\check x_{m}),
$$
where $\check x_{i}=(x_{l_{i-1}+1},\ldots, x_{l_{i}})$
and set
$$
\check \Omega^{i}=\Omega^{l_{i-1}+1}
 \times\cdots\times \Omega^{l_{i}},\quad
\hat \Omega^{i}=
\Omega^{l_{i-1}+1}\times\cdots\times \Omega^{d},
$$
 $\hat x_{i}=(x_{l_{i}+1},\ldots,x_{d})$.
Take
$k(1),\ldots,k(d)\in\{1,2,\ldots\}$ and, for $n\in\bZ$,
let
$$
 \check C_{n}^{i}=[0,2^{-nk(l_{i-1}+1)})\times\cdots\times
[0,2^{-nk(l_{i})})
$$
 be a subset
of $\check \Omega^{i}$ and $C_{n}=\check C_{n}^{1}\times\cdots\times
\check C_{n}^{m}$. By $A_{p}$-weights on $\check \Omega^{i}$
we mean the $A_{p}$-weights relative to all translates
of $\check C_{n}^{i}$, $n\in\bZ$, belonging to $\check \Omega^{i}$, and, naturally,
$A_{p}$-weights on $\Omega $ are defined using all
translates of $C_{n}$, $n\in\bZ$, belonging to $ \Omega$.

 \begin{theorem}
                           \label{theorem 12.2.1}
Let $K_0,p_{k}\in(1,\infty)$, $w^k\in
A_{p_{k}}( \check\Omega^k)$,
 $[w^k]_{p_{k}}\le K_0$, $k=1,\ldots,m$, and
$u,g$ be measurable functions on $\Omega$. Then there exists
a constant $\Lambda_0=\Lambda_0(d,p_1,\ldots,p_{m},
k(1),\ldots,k(d),K_0)\ge
1$ such that if
$$
\|u\|_{L_{p_{1}}(w\,d\mu)}
\le N_0\|g\|_{L_{p_{1}}(w\,d\mu)}
$$ for some $N_0\in (0,\infty)$ and for every $w\in
A_{p_{1}}(\Omega)$ with
 $[w]_{p_{1}}\le \Lambda_0$, then we have
$$
\|u\|_{L_{p_{1},\ldots,p_{m}}(w^1,\ldots,w^{m})}\le
N\|g\|_{L_{p_{1},\ldots,p_{m}}(w^1,\ldots,w^m)},
$$ where  the norms are defined as in \eqref{4.9.6}
replacing $dx_{i}$ by $w^{i}(\check x_{i})\,d\check x_{i}$, the constant
$N$ depends only on $d$, $p_1,\ldots,p_m,
k(1),\ldots,k(d)$, $K_0$, and
$N_0$.
\end{theorem}

Proof. We follow the proof of Corollary 2.7
in \cite{DK16}.
Recall
the extrapolation theorem of
 J. L. Rubio de Francia
\cite{MR745140} which says that  for any constant
$\Lambda_{j }\in (1,
\infty)$, $j= 1,\ldots,m$, there exists a constant
$\Lambda_{j-1}=\Lambda_{j-1}(d-j,p_{j},p_{j+1},
 K_{0} \Lambda_{j } )\in ( 1,
\infty)$ (we drop its dependence on the $k(i)$'s)
such that, if

(a) for two nonnegative functions
$U_{j}$ and $G_{j}$ on $\hat\Omega^{j+1}$ it holds that
\begin{equation}
                                         \label{4.22.3}
\int_{\hat\Omega^{j+1} }
U_{j}^{p_{j}}w(\hat x_{j+1} )\,
d\hat x_{j+1}   \leq N_{j}\int_{\hat\Omega^{j+1} }
G_{j}^{p_{j}}w(\hat x_{j+1} )\,
d\hat x_{j+1}
\end{equation}
 for some $N_j\in (0,\infty)$ and for every $w\in
A_{p_{j}}(\hat\Omega^{j+1} )$ with
 $[w]_{p_{j}}\le \Lambda_{j-1}$, then

(b) we have
\begin{equation}
                                         \label{4.22.4}
\int_{\hat\Omega^{j+1} }
U_{j}^{p_{j+1}}w(\hat x_{j+1} )\,
d\hat x_{j+1}   \leq N_{j+1}\int_{\hat\Omega^{j+1} }
G_{j}^{p_{j+1}}w(\hat x_{j+1} )\,
d\hat x_{j+1}
\end{equation}
 for some $N_{j+1}\in (0,\infty)$,
depending only on $d$, $j$,   $K_0\Lambda_j$, $p_{j}$, $p_{j+1}$,
and $N_{j}$,
 and for every $w\in
A_{p_{j+1}}(\hat\Omega^{j+1} )$ with
 $[w]_{p_{j+1}}\le  K_{0}\Lambda_{j}$.

In this form the theorem
is proved in \cite{DK16}.
We define $\Lambda_{m-1}=1$ and find all $\Lambda_{j}$,
$j=0,1,\ldots,m-1$.  Then assume that $m\geq 2$ and
define   $U_{0}(x)=u(x)$,
$$
U_{j}(\hat x_{j+1})=
\Big(\int_{\check \Omega^{j}}U^{p_{j}}_{j-1}(\hat x_{j})
\, w^{j}(\check x_{j})\,d\check x_{j}\Big)^{1/p_{j}},\quad  1\leq j\leq m-1,
$$
and similarly we introduce
$G_{j}$'s
by taking $g$ in place of $u$. To prove the theorem,
it suffices to prove that (b) holds for
$j=m-1$ because $w^{m}\in A_{p_{m}}(\check \Omega^{m})$
and $[w^{ m}]_{p_{m}}\leq K_{0}
= K_{0}\Lambda_{m-1}$. We are going to use the induction on $j=0,1,\ldots,m-1$.

Observe that
 (b) holds for $j=0$   by assumption.
Suppose that it holds for a $j\in \{0,1,\ldots,m-2\}$.
Then \eqref{4.22.4} also holds for
$$
w(\hat x_{j+1} ):=w^{j+1}(\check x_{j+1})w(\hat x_{j+2} )
$$
if $w^{j+1}\in A_{p_{j+1}}( \check\Omega^{j+1})$
and $w(\hat x_{j+2} ) \in A_{p_{j+1}}(\hat \Omega^{j+2})$ with
$$
[w^{j+1}]_{p_{j+1}}\le K_{0},\quad
[w(\hat x_{j+2} )]_{p_{j+1}}\le \Lambda_{j }
$$
because then
$
[w(\hat x_{j+1} )]_{p_{j+1}}\le K_{0}\Lambda_{j }
 $. Remarkably, this implies that
(a) holds with $j+1$ in place of $j$.
Then (b) also holds with $j+1$ in place of $j$.
This justifies the induction and proves the theorem.

%\bibliographystyle{plain}
%\bibliography{fullynonlinear}

\begin{thebibliography}{10}

\bibitem{BP61}
A.~Benedek and R.~Panzone.
\newblock The space {$L\sp{p}$}, with mixed norm.
\newblock {\em Duke Math. J.}, 28:301--324, 1961.

\bibitem{CC_95}
Luis~A. Caffarelli and Xavier Cabr\'e.
\newblock {\em Fully nonlinear elliptic equations}, volume~43 of {\em American
  Mathematical Society Colloquium Publications}.
\newblock American Mathematical Society, Providence, RI, 1995.

\bibitem{CD}
Mar\'\i a~E. Cejas and Ricardo~G. Dur\'an.
\newblock Weighted a priori estimates for elliptic equations.
\newblock {\em Studia Math.}, 243(1):13--24, 2018.

\bibitem{CKS00}
M.~G. Crandall, M.~Kocan, and A.~\'Swi{\c e}ch.
\newblock {$L^p$}-theory for fully nonlinear uniformly parabolic equations.
\newblock {\em Comm. Partial Differential Equations}, 25(11-12):1997--2053,
  2000.

\bibitem{DK16}
Hongjie Dong and Doyoon Kim.
\newblock On {$L_p$}-estimates for elliptic and parabolic equations with
  {$A_p$} weights.
\newblock {\em  Trans. Amer. Math. Soc.}, 370(7):5081--5130, 2018.

\bibitem{DG18}
Hongjie Dong and Chiara Gallarati.
\newblock Higher-order parabolic equations with vmo assumptions and general
  boundary conditions with variable leading coefficients.
\newblock {\em International Mathematics Research Notices}, page rny084, 2018.



\bibitem{DKL_12}
Hongjie Dong, N.~V. Krylov, and Xu~Li.
\newblock On fully nonlinear elliptic and parabolic equations with {VMO}
  coefficients in domains.
\newblock {\em Algebra i Analiz}, 24(1):53--94, 2012.

\bibitem{Gr09}
Loukas Grafakos.
\newblock {\em Modern {F}ourier analysis}, volume 250 of {\em Graduate Texts in
  Mathematics}.
\newblock Springer, New York, second edition, 2009.

\bibitem{NK_09}
Vladimir Kozlov and Alexander Nazarov.
\newblock The {D}irichlet problem for non-divergence parabolic equations with
  discontinuous in time coefficients.
\newblock {\em Math. Nachr.}, 282(9):1220--1241, 2009.

\bibitem{Kr_01}
N.~V. Krylov.
\newblock The heat equation in {$L_q((0,T),L_p)$}-spaces with weights.
\newblock {\em SIAM J. Math. Anal.}, 32(5):1117--1141, 2001.

\bibitem{Kr_08}
N.~V. Krylov.
\newblock {\em Lectures on elliptic and parabolic equations in {S}obolev
  spaces}, volume~96 of {\em Graduate Studies in Mathematics}.
\newblock American Mathematical Society, Providence, RI, 2008.

\bibitem{Kr_13}
N.~V. Krylov.
\newblock On the existence of {$W_p^2$} solutions for fully nonlinear elliptic
  equations under relaxed convexity assumptions.
\newblock {\em Comm. Partial Differential Equations}, 38(4):687--710, 2013.

\bibitem{Kr_08_1}
N.~V. Krylov.
\newblock On parabolic {PDE}s and {SPDE}s in {S}obolev spaces {$W^2_P$} without
  and with weights.
\newblock In {\em Topics in stochastic analysis and nonparametric estimation},
  volume 145 of {\em IMA Vol. Math. Appl.}, pages 151--197. Springer, New York,
  2008.

\bibitem{Kr_10}
N.~V. Krylov.
\newblock On {B}ellman's equations with {VMO} coefficients.
\newblock {\em Methods Appl. Anal.}, 17(1):105--121, 2010.

\bibitem{Kr_18}
N.~V. Krylov.
\newblock On the existence of $W^{1,2}_{p}$
solutions for fully nonlinear parabolic equations under either relaxed
or no  convexity assumptions.
\newblock{accepted for CMSA Nonlinear Equation
Publication, arXiv:1705.02400}

\bibitem{Kr17}
N.~V. Krylov.
\newblock{\em Sobolev and viscosity solutions
for fully nonlinear  elliptic and parabolic
 equations}.
\newblock to appear with AMS.


\bibitem{Lin86}
Fang-Hua Lin.
\newblock Second derivative {$L^p$}-estimates for elliptic equations of
  nondivergent type.
\newblock {\em Proc. Amer. Math. Soc.}, 96(3):447--451, 1986.

\bibitem{MR745140}
Jos\'e~L. Rubio~de Francia.
\newblock Factorization theory and {$A\sb{p}$}\ weights.
\newblock {\em Amer. J. Math.}, 106(3):533--547, 1984.

\bibitem{Wi09}
Niki Winter.
\newblock {$W^{2,p}$} and {$W^{1,p}$}-estimates at the boundary for solutions
  of fully nonlinear, uniformly elliptic equations.
\newblock {\em Z. Anal. Anwend.}, 28(2):129--164, 2009.

\end{thebibliography}

\end{document}